\RequirePackage{fix-cm}
\documentclass[smallextended,natbib]{svjour3}       
\smartqed  
\usepackage{graphicx}
\usepackage{mathptmx}      
\usepackage{amssymb}
\usepackage{amsmath}
\usepackage{latexsym}
\usepackage{subscript}
\usepackage{mathrsfs}
\usepackage[parfill]{parskip}
\usepackage[all]{xy}

\usepackage[textwidth=1in,textsize=footnotesize]{todonotes}

\newcommand*{\mP}{\mathfrak{P}}
\newcommand*{\mQ}{\mathfrak{Q}}
\newcommand*{\mR}{\mathfrak{R}}
\newcommand*{\mS}{\mathfrak{S}}
\newcommand*{\mM}{\mathfrak{M}}
\newcommand*{\Lang}{\mathcal{L}}
\newcommand*{\Iso}{\mathrm{Iso}}
\newcommand*{\For}{\mathrm{For}}

\journalname{Synth\`ese}

\begin{document}

\title{Carnap's Early Metatheory: Scope and Limits}

\author{Georg Schiemer \and Richard
  Zach  \and Erich Reck}

\institute{Georg Schiemer \at
  Department of Philosophy, 
  University of Vienna, Universit\"atsstra\ss e 7, 1010 Vienna, Austria \& 
\at
Munich Center for Mathematical Philosophy, LMU Munich,
Geschwister-Scholl-Platz 1, 
D-80539 \\ Munich, Germany \\
  \email{georg.schiemer@univie.ac.at}
  \and
  Richard Zach \at
  Department of Philosophy, 
  University of Calgary, 2500 University Drive NW, Calgary AB T2N 1N4, Canada\\
  \email{rzach@ucalgary.ca}
  \and
  Erich Reck \at Department of Philosophy, 
  University of California Riverside, 900 University Avenue,
Riverside, CA 92521, USA \\
  \email{erich.reck@ucr.edu}
}

\date{Version: \today  / Received: date / Accepted: date}

\maketitle

\begin{abstract}
In his \emph{Untersuchungen zur allgemeinen Axiomatik} (1928) and
\emph{Abriss der Logistik} (1929), Rudolf Carnap attempted to formulate 
the metatheory of axiomatic theories within a single, fully interpreted 
type-theoretic framework and to investigate a number of meta-logical 
notions in it, such as those of model, consequence, consistency, completeness, 
and decidability.  These attempts were largely unsuccessful, also in his
own considered judgment. A detailed assessment of Carnap's attempt shows, 
nevertheless, that his approach is much less confused and hopeless than 
it has often been made out to be.  By providing such a reassessment, the 
paper contributes to a reevaluation of Carnap's contributions to the 
development of modern logic.
\end{abstract}

\section{Introduction}

Rudolf Carnap's contributions to logic prior to his \emph{Logische
Syntax der Sprache} \citeyearpar{Carnap1934} consist, among other
things, of work on the formalization of mathematical theories and
their metamathematical investigation. His main contributions to this
topic are contained in his logic textbook \emph{Abriss
der Logistik} \citeyearpar{Carnap1929} and the manuscript
\emph{Untersuchungen zur Allgemeinen Axiomatik}, 
written in Vienna around 1928, but unpublished until the edition
\citep{Carnap2000}.\footnote{Carnap presented the main
results of his \emph{Untersuchungen} manuscript at the \emph{First
Conference on Epistemology of the Exact Sciences}, August 1929,
in Prague. Already at the end of 1928/beginning of 1929, he
considered submitting a version of the manuscript for publication,
and circulated it to mathematicians and logicians, including Baer,
Behmann, G\"odel, H\"arlen, and Fraenkel; see
\cite{AwodeyCarus2001}, note~1.  Behmann provided extensive
comments; some of his marginalia are preserved in one of the
versions of the typescript in Carnap's papers (ASP RC 080--34--02),
and longer comments in Behmann's papers (Staatsbibliothek zu Berlin,
Nachla\ss{} 335 (Behmann), K.~1 I~10). Behmann at first suggested the manuscript
was not ready to be published in a pure mathematics journal, but
withdrew these objections after closer reading in March 1929.
However, Carnap delayed revision of the manuscript and eventually,
in 1930, abandoned the plan to publish it, essentially after some 
conversations with Tarski during his first visit to Vienna; 
see~\citet{Reck2007}. The edition \citep{Carnap2000} includes items RC
080--34--02, 080--34--03, and 080--34--04 from Carnap's papers, with
uniformized page, section, and theorem numbering.  These original
items are now available online from the Archives of Scientific
Philosophy (www.library.pitt.edu/rudolf-carnap-papers). We 
provide references to both the (out-of-print) edition
and the original typescripts.  Carnap's work on general axiomatics
is also documented in \cite{Carnap1927b,Carnap1930} and
\cite{CarnapBachmann1936}. For a general overview of logical work on
axiom systems before Carnap, see \cite{MancosuZachBadesa2009}.}
Carnap's study of the metatheory of axiomatic theories developed in
these texts remains interesting today. On the one hand, it is highly
original, especially given the fact that his results were formulated
prior to G\"odel's incompleteness results and to Tarski's work on
truth and logical consequence from the mid-1930s. It thus illustrates
Carnap's still often neglected contributions to the development 
of modern logic.  On the other hand, his approach has interesting 
points of contact with parallel contributions to logic from the same 
period. This concerns, in particular, Tarski's earlier work on the
``methodology of the deductive sciences'', from the 1920s and early
1930s. Like Tarski, Carnap aims to give an explication of several 
metatheoretical concepts that play a role in modern axiomatics
and to specify their logical relations.

The general aim of the present paper is to reassess 
Carnap's resulting proposals in \emph{Untersuchungen} and
\emph{Abriss}.\footnote{Carnap's more general views on logic were in
flux during the 1920s and early 30s, especially after he encountered 
\cite{Hilbert1928}, \cite{HilbertAckermann1928}, \cite{Godel1929}, 
and Tarski's early meta-logical work.  Our discussion focuses on 
\cite{Carnap2000} and \cite{Carnap1929}, two texts that were
composed largely before those influences became dominant.  It should
be noted, however, that Carnap followed the developments of logic in
Hilbert's school closely even before the appearance of Hilbert and
Ackermann's textbook.  In fact, the \emph{Abriss} was originally
conceived as a joint project with Hilbert's student Heinrich Behmann
in the early 1920s.  However, Behmann came to prefer his own, more
algebraic notation (see \citealt{MancosuZach2015}), while Carnap
thought the Peano-Russell-Whitehead notation made popular by
\emph{Principia} was preferable for a textbook presentation (Carnap
to Behmann, February 19, 1924, Staatsbibliothek zu Berlin,
Nachla\ss{} 335 (Behmann), K.~1 I~10). For additional historical
background, see \cite{AwodeyCarus2001}, \cite{AwodeyReck2002},
and \cite{Reck2004,Reck2007}.  A broader study of the evolution 
of Carnap's metatheory up to and beyond the 1930s
cannot be undertaken here; it deserves 
a separate paper.}  To do this, we articulate---in more detail than the 
literature on this topic has done so far---the logical machinery of and 
the conceptual assumptions behind Carnap's attempt to formulate the 
metatheory of axiomatic theories in a type-theoretic framework. This 
involves addressing questions such as the following: How precisely 
did Carnap explicate metalogical notions in \emph{Untersuchungen}, 
such as those of model, logical consequence, and completeness? In 
which ways are the results related to similar contributions by his 
contemporaries Behmann, Fraenkel, Hilbert, and Tarski, among others? 
In what ways does his account differ most significantly from today's 
metalogic? And are there any important structural similarities between 
them, despite the conceptual differences?

In pursuing this objective, we focus on three characteristic 
features of Carnap's approach in \emph{Untersuchungen} that 
distinguish it from metalogic as practiced today, at least at first
glance. The first feature concerns the fact that a clear-cut
distinction between semantics and proof theory is still missing in
Carnap's text, as becomes evident at various points. The second
feature to be reconsidered is his particular type-theoretic approach:
Carnap attempts to use the same logical language both for formulating 
axiom systems and their consequences and for formulating metatheoretic
concepts and metatheorems about such systems. Because of this feature,
Carnap has repeatedly been accused of failing to make a necessary 
distinction between object language and metalanguage, an accusation 
we will reassess. A third feature characteristic of Carnap's 
pre-\emph{Syntax} logic concerns his treatment of logical 
languages themselves, which are not conceived of as formal
or disinterpreted by him, as is usual today, but as ``meaningful
formalisms'' that come with a fixed and intended interpretation.

In our reevaluation of Carnap's early metatheory we aim to avoid 
simple, anachronistic \emph{ex post} corrections of the logical 
``mistakes'' contained in Carnap's work, i.e., to point to differences 
between his and the now standard approach only to dismiss
Carnap's. Taking Carnap's original definitions seriously, analyzing
them in their historical context, and evaluating them, both on their 
own terms and in comparison to later results, is meant to contribute
to a substantive, balanced, and genuinely contextual history 
with respect to Carnap's contributions to logic. More particularly, 
our aim is to expand upon the existing literature by providing 
detailed treatments of Carnap's notions of logical consequence,
model, isomorphism, formality, and various ``k-concepts.''  But this 
will have repercussions on other issues as well, including the
notions of completeness and decidability.\footnote{Our
discussion in this paper builds on existing scholarship. The first
phase of engagement with Carnap's ``general axiomatics'' project
consisted in the pioneering but overly critical work by
\cite{Coffa1991} and \cite{Hintikka1991,Hintikka1992}. Both authors
criticize the ``monolinguistic approach'' adopted by Carnap, i.e.,
his attempt to express axiomatic theories and their metatheory in a
single type-theoretic language. A second phase set in with Awodey's
and Carus' important paper on the logical and philosophical analysis
of the main technical result in \emph{Untersuchungen}, the so-called
\emph{Gabelbarkeitssatz} \citep{AwodeyCarus2001}. That paper was
followed by a number of articles aimed at a more balanced account of
Carnap's project. In them, not only the limitations of his approach
were acknowledged, but also its innovative aspects and its
significant influence on later developments in metalogic
\citep{AwodeyReck2002,Reck2004,Goldfarb2005,Reck2007}. The third
phase consists in fairly recent scholarship in which attention is
drawn to previously neglected details of Carnap's early model theory
and in which its role in the development of metalogic is spelled out more
\citep{Reck2011,Schiemer2012,Schiemer2012b,Schiemer2013,SchiemerReck2013,Loeb2014,Loeb2014b}.}

The paper is organized as follows. We begin, in section~\ref{Axiomatic}, 
by introducing the approach to axiomatic theories and their 
models as developed by Carnap in the \emph{Abriss}, which 
also provides the framework for his metalogical investigations in
\emph{Untersuchungen}. We then introduce the main metatheoretic 
notions as defined in \emph{Untersuchungen}, focusing especially 
on Carnap's explication of logical consequence and (various versions 
of) completeness.\footnote{The content of this section will be familiar 
to the Carnap specialist. But it can serve as an introduction to 
the reader not yet familiar with Carnap's general axiomatics program 
and the literature around it; we also expand on that literature with 
respect to some details. Similar remarks apply to the next two sections.}  
In section~\ref{Completeness}, we discuss in what way Carnap's approach
to semantics in \emph{Untersuchungen} differs from the now standard
approach.  We aim to show that Carnap's approach was neither
idiosyncratic in its historical context, nor is it as misguided as
was often claimed later.  In section~\ref{Domain} we address 
a specific question regarding Carnap's semantic notions, that of
domain variation. This is a crucial aspect of the current definitions of
model and consequence and, contrary to what has often been
assumed, it can be accommodated within Carnap's type-theoretic framework.  
Next, in section~\ref{Consequence}, we reconsider Carnap's notion
of logical consequence in relation to the proof-theoretic notion of
derivability. It is here that Carnap's approach is hardest to reconcile 
with the now standard perspective, and we attempt to diagnose the
resulting difficulties precisely. In section~\ref{Isomorphism}, we
focus on the notion of higher-level isomorphism, the closely related
notion of ``formality'', as part of his definitions of forkability 
and completeness, and Carnap's proof of the
\emph{Gabelbarkeitssatz}.  In section~\ref{akConcepts}, we further
elaborate on the lack a proper distinction between model-theoretic and
proof-theoretic notions in \emph{Untersuchungen}, here especially the 
difficulties caused for Carnap's discussion of the distinction 
between ``absolute forkability'' and ``constructive decidability'', or
between ``absolute'' and ``constructive'' concepts more generally.
Finally, Section \ref{Summary} contains a brief summary of our
findings, together with suggestions for possible future research.

\section{Axiomatic theories and Carnap's model theory}
\label{Axiomatic}

Carnap's approach to formalizing axiomatic theories is discussed in
detail in \emph{Abriss} and \emph{Untersuchungen}. In this section, 
we give a brief and slightly updated presentation of this basic 
approach. According to Carnap, axioms, axiom systems, and 
corresponding theorems are to be expressed in a higher-order language,
more precisely a language of simple (or de-ramified) types.\footnote{A
detailed presentation of Carnap's type-theoretic logic can be found 
in \citep[\S 9]{Carnap1929}.} A formal axiom system is understood by
him as a ``theory schema, the empty form of a possible theory''. Here
the primitive terms of a theory are represented by means of variables
(of given arity and type) $X_1, \dots, X_n$ (and not, as is usual
today, by means of schematic nonlogical constants). Axioms, axiom 
systems, and their theorems are then symbolized as propositional
functions of the form $f(X_1, \dots, X_n)$, i.e. as open formulas in 
the current sense. Finally, a corresponding theory is a formula that 
consists of the conjunction of the relevant axioms (and not, as is
usual today, of the closure of the set of axioms under deductive or
semantic consequence).\footnote{To be more precise, Carnap 
describes two different ways of formalizing axiom systems in \emph{Abriss}. 
According to the first approach, primitive terms are expressed in terms of 
nonlogical constants that already have a predetermined (and fixed) meaning. 
According to the second approach (the one typical for formal axiomatics), 
the primitive terms defined by a theory have no predetermined meaning; 
they express, in Carnap's terminology, only ``improper concepts''. Such 
terms are symbolized by means of variables in the way described above. 
Compare \cite[\S 30b]{Carnap1929}. In what follows, we will 
always assume this second approach to the formalization of theories.} 

We can consider one of Carnap's own examples to further illustrate this 
approach. In the second part of \emph{Abriss}, axiom systems from
arithmetic, geometry, topology, and physics are formalized in a
type-theoretic language. Among them are Hausdorff's axioms for
topology and the corresponding topological theory. The single ``primitive 
sign'' of this theory is a binary relation variable $\mathsf{U}(\alpha, x)$ 
that stands for ``$\alpha$ is a neighborhood of $x$''. A second set 
variable $\mathsf{pu}$ (standing for the set of points) is then
defined relative to $\mathsf{U}(\alpha, x)$, namely as the range 
of the neighborhood relation.  The theory's axioms are the 
following:\footnote{Put informally, these axioms state that (1a)
neighborhoods are classes of points, (1b) a point belongs to each of
its neighborhoods; (2) the intersection of two neighborhoods 
of a point contains a neighborhood, (3) for every point of a
neighborhood~$\alpha$ a subclass of~$\alpha$ is also a neighborhood,
and (4) for two different points there are two corresponding
neighborhoods that do not intersect; see \citet[\S
33a]{Carnap1929}. Carnap's own type-theoretic presentation of
these axioms has been slightly modernized here, although we 
have also preserved some of his idiosyncratic notation.}
\begin{align}
\tag{Ax1a}  & \mathrm{Dom}(\mathsf{U}) \subset \wp(\mathsf{pu})\\
\tag{Ax1b} & \mathsf{U}  \subset \ \in^{-1} \\
\tag{Ax2} &\forall \alpha, \beta, x \ ( \mathsf{U}(\alpha,x)  \land 
\mathsf{U}(\beta, x)  \rightarrow \exists \gamma ( \mathsf{U}(\gamma,x) 
\land \gamma \subset (\alpha \cap \beta)))\\
\tag{Ax3} &\forall \alpha , x \ ( \alpha \in \mathrm{Dom}(\mathsf{U}) 
 \land x \in \alpha \rightarrow \exists \gamma ( \mathsf{U}(\gamma,x) 
\land \gamma \subset \alpha))\\
\tag{Ax4} &\forall x, y (x, y \in \mathsf{pu} \land x  \neq y \rightarrow \exists \alpha, \beta ( \mathsf{U}(\alpha,x) \land 
\mathsf{U}(\beta,y)  \land (\alpha \cap \beta) = \emptyset))
\end{align}
Notice that each of these axioms is purely logical, i.e. an open
formula containing only logical primitives of the type-theoretic
language and logico-mathematical constants (such as `$\mathrm{Dom}(X)$' or `$x
\in \alpha$') that are explicitly defined in that language. In
addition, the axioms contain two highlighted ``primitive signs,''
i.e., free variables that represent the primitive terms of the
theory. Finally, Hausdorff's theory is expressed by an open formula
$\Phi_{Hausd}$ that consists of the conjunction of
$Ax1-Ax4$.

In \emph{Abriss}, Carnap is not yet explicit about how he understands
interpretations or models for such theories. However, in
\emph{Untersuchungen}, as well as in \citet{Carnap1930} and
\citet{CarnapBachmann1936}), a detailed discussion of such models is
provided. Roughly, models of a theory are treated as ordered sequences
of predicates of the type-theoretic language that can be substituted
for the variables representing the ``primitive terms'' of a
theory. (Model classes, i.e., the classes of all models corresponding
to a given theory, can then be considered correspondingly.) Compare,
for instance, Carnap's following characterization:
\begin{quote}
If $f\mR$ is satisfied by the constant $\mR_1$, where $\mR_1$ is an
abbreviation for a system of relations $P_1$, $Q_1$, \dots, then
$\mR_1$ is called a ``\emph{model}'' of~$f$. A model is a system of
concepts of the basic discipline, in most cases a system of numbers
(number sets, relations, and such). \citep[p.~303]{Carnap1930}
\end{quote}
It is not fully clear from this passage (nor from related ones) whether 
Carnap understands models as linguistic entities or, instead, as 
set-theoretic entities, i.e., relational structures in the current sense 
(which could be specified in his type-theoretic background language 
by using constants and predicates of the right type). The important 
point to note, in any case, is that satisfaction---truth in a model---is 
treated substitutionally by Carnap: Given a sequence of predicates 
of the form $\mM = \langle R_1, \dots, R_n \rangle$, where each 
predicate $R_i$ represents an admissible value of the $i$-th primitive 
sign of the theory $f(X_1, \dots, X_n)$, he maintains that the theory 
is true in the model $\mathfrak{M}$ if the sentence resulting from 
the substitution of each primitive variable by its corresponding 
predicate in the model, i.e., $f(R_1, \dots, R_n)$, ``is true'', or ``holds'', 
in his underlying type-theoretic logic. (The central interpretive issue 
of what such type-theoretic ``truth'' or ``holding'' amounts to for 
Carnap will be discussed in the next section.)

Returning to the above example, we can say that a particular model 
of Hausdorff's topology is given by a tuple $\langle \mathsf{pu},
\mathsf{U} \rangle$ consisting of a unary and a binary predicate
that---if substituted for the two primitive variables---turn axioms
Ax1--Ax4 into true sentences. If we put aside Carnap's substitutional
approach to truth for the moment, his account of formal axiomatic 
theories and their models looks, on the surface, like the now
standard model-theoretic treatment. However, additional and 
important differences become evident if one pays closer attention 
to how exactly metatheoretic concepts are defined by Carnap 
within this type-theoretic framework.

Carnap's main contribution in \emph{Untersuchungen} is not 
the logical formalization of axiomatic theories presented in the 
previous section. Carnap's way of formalizing theories in terms of
propositional functions, and of their primitives in terms of free
variables, was common at the time; similar approaches can be
found in work by logicians like Russell, Tarski, and Langford, and
mathematicians like Peano, Hilbert, and Huntington. Rather, it consists 
of his explications of several metatheoretic notions that had been used 
informally in mathematics since Dedekind and Hilbert.\footnote{See 
\cite{AwodeyCarus2001}, \cite{Reck2004}, and \cite{Carus2007} 
for the broader philosophical context of Carnap's metatheoretical work 
in \emph{Untersuchungen}, including its connections to the \emph{Aufbau}
\citeyearpar{Carnap1928a}. Carnap's own views on the broader philosophical significance 
of his \emph{general axiomatics} project are elaborated, in particular, 
in ``Eigentliche und uneigentliche Begriffe'' \citep{Carnap1927b}.  
Further exploring this part of the background deserves a separate
investigation, one that may help to account more for some of the 
idiosyncracies in Carnap's metalogical approach even beyond 1930.
However, we cannot take on this 
additional task in the present paper.}  A closer look at these explications reveals 
several striking differences to the now standard approach.  Most importantly, 
in contrast to current model theory and proof theory, Carnap provided
metatheoretic definitions for the theories he was considering partly
in the same language in which they were formulated, and partly in an
informal metalanguage (German). He did not clearly distinguish between
an object language (in which theories are expressed) and a
metalanguage (in which model-theoretical and proof-theoretical results
for such theories are stated). To a certain degree, it is justified to
say that Carnap's project is ``monolinguistic'' in character: it
assumes a single higher-order language, more specifically a language
of simple type theory (henceforth $\Lang_{TT}$), in which both
axiomatic theories and their metatheoretical properties are
expressed.\footnote{It should be emphasized, however, that the early 
critical discussion of this monolinguistic character by Coffa and Hintikka
does not do justice to the subtleties of Carnap's approach. As was
first analyzed in detail in \cite{AwodeyCarus2001}, Carnap
presupposes a ``basic system'' (\emph{Grunddisziplin}), that is, a
theory of ``contentful'' logical and mathematical concepts, that
forms the background for the logical study of axiomatic theories. In
\cite{Carnap2000}, simple type theory is proposed as one possible
choice for such an interpreted basic system. As we will show
below, Carnap employed it similarly to the way in which we
use axiomatic set theory as the background for the metatheoretic
study of axiomatic theories today.  Thus, while a clear-cut
object/metalanguage is still missing in \emph{Untersuchungen},
Carnap's approach can be characterized as metatheoretic.  Compare
\cite{SchiemerReck2013} for a further discussion of this point. In 
addition, it should be stressed that Carnap's convention to formulate both
axiom systems and their metatheoretic properties within a single and
interpreted \emph{Grunddisziplin} was not uncommon at the time. 
See, for instance, a related discussion of a necessary ``absolute
foundation'' of axiomatic theories in Fraenkel's third edition of
\emph{Introduction to Set Theory} \citep[p.154]{AwodeyCarus2001}.  
Finally, a similar idea is present in own Carnap's subsequent work, 
in particular in his \emph{Logical Syntax}: it is shown here that the
``syntax language'' for language LI can be formulated within LI
itself. See \citep[\S 18]{Carnap1934}.}

Let us now consider some particular metatheoretic notions discussed in
\emph{Untersuchungen} so as to illustrate Carnap's approach further. 
The notion of ``consequence'' is introduced as follows by him:
\begin{definition}[Consequence]
A sentence $g$ is a consequence of axiom system $f$ iff
\[
\forall \mR(f  \mR \rightarrow
g \mR)
\]
\end{definition}
Consequence is here specified in terms of the material conditional, or
more precisely, via a quantified conditional statement
of~$\Lang_{TT}$.\footnote{In Carnap's own words: ``$g$ is called a
  `\emph{consequence}' of $f$, if $f$ generally implies $g$: $\forall
  \mR(f\mR \rightarrow g\mR$), abbreviated: $f \rightarrow g$. The
  consequence is, as is the AS, not a sentence, but a propositional
  function; only the associated implication $f \rightarrow g$ is a
  sentence, namely a purely logical sentence, thus a tautology, since
  no nonlogical constants occur.''  \citep[p.~304]{Carnap1930}} This
differs obviously from the now standard model-theoretic approach to 
the notion; but there are also structural similarities. Notice, in particular,
that Carnap's definition involves universal quantification over what
he takes to be models of a theory. In addition, in \citep{Carnap2000}
he states explicitly that his notion of consequence is not to be
conflated with Hilbert's notion of derivability in a formal
system.\footnote{In \emph{Untersuchungen} Carnap argues that ``$g$
  follows from $f$'' and ``$g$ is derivable from $f$ in TT'', while
  not identical, are equivalent. See \citep[p.~92]{Carnap2000}; ASP RC
  080--34--03, p.~41a--b. We will come back to this issue later.}

Based on this definition of consequence, together with related
type-theoretic definitions of further auxiliary notions (such as those
of isomorphism between models and satisfiability of a theory,
more on which below), Carnap specifies three notions of completeness 
for a theory, namely ``monomorphicity'', ``(non-)forkability'', and 
``decidability''.\footnote{As discussed by \cite{AwodeyCarus2001} 
and \cite{AwodeyReck2002}, these correspond roughly to what would 
today be called categoricity, semantic completeness, and syntactic
completeness. Then again, this correspondence must be viewed with
caution, especially in the case of the latter two notions, as we will 
shown in section~\ref{akConcepts}.} The three notions were 
identified earlier by \cite{Fraenkel1928}.  Carnap is here responding 
to an open research question posed by Fraenkel, viz., what the relationship 
between them is.  His specification of monomorphicity is the following:
\begin{definition}[Monomorphic]
An axiom system $f$ is \emph{monomorphic} iff
\[
\exists \mR (f \mR)  \land  \forall \mP \forall \mQ [(f\mP \land f\mQ) 
\rightarrow Ism_{q}(\mP, \mQ)]
\]
\end{definition}
Here the first conjunct states, along Carnapian lines, that the theory
$f$ is satisfied, while the second conjunct asserts that all of its
models are isomorphic.

A second version of completeness of an axiom system is that of
``non-forkability'' (``\emph{Nicht-Gabelbarkeit}''). With respect to it,
the crucial Carnapian definition is the following:
\begin{definition}[Forkability]
An axiom system $f$ is \emph{forkable at a sentence~$g$} iff
\[
\exists \mR (f \mR \land g \mR) \land \exists \mS (f \mS \land \lnot
g\mS) \land \underbrace{\forall \mP \forall \mQ [(g\mP \land
    Ism_{q}(\mP, \mQ)) \rightarrow g\mQ]}_{\For(g)}
\]
\end{definition}
The first two conjuncts in this formula state that theory $f$ is
satisfiable jointly with sentence (or propositional function)~$g$ as
well as with its negation~$\lnot g$. The third conjunct expresses the
fact that $g$ is ``formal'': if $g$ is true in a particular model,
then it is also true in any other model isomorphic to the former. (We
will come back to this notion of formality below.)  A theory is then
``non-forkable'' if it is not forkable in this sense. In other words,
it is non-forkable if there is no formal sentence~$g$ such that both
$f \land g$ and $f \land \lnot g$ are satisfiable.

Carnap's definition of ``decidability'' (``\emph{Entscheidungsdefinitheit}'') 
is this: 
\begin{definition}[Decidability]
An axiom system $f$ is \emph{decidable} iff
\[
\exists \mR (f \mR)  \land  \forall g (For(g) \rightarrow \forall \mP 
((f\mP \rightarrow g\mP) \lor (f\mP \rightarrow \lnot g\mP)))
\]
\end{definition}
This definition states that a theory is decidable if it satisfiable and if 
for any formal sentence $g$ (again expressed by a propositional
function like above) either it or its negation is a consequence of the
theory. Note that, given Carnap's definition of consequence, this 
makes it another version of the later notion of semantic completeness, 
rather than of that of decidability in the sense of computability theory. 
(For a closely related discussion of Carnap's ``constructive'' variant 
of this notion, i.e. ``k-decidability'', see section~\ref{akConcepts}.)

Monomorphicity, (non-)forkability, and decidability are the three
central notions discussed in the first part of \emph{Untersuchungen} 
\citep{Carnap2000}. Even from today's vantage point, these notions 
are genuinely metatheoretic: Each of them captures a different 
logical property of axiomatic theories in terms of their consequences 
and their models.\footnote{The central metatheorem in \citep{Carnap2000}, 
Carnap's so-called \emph{Gabelbarkeitssatz}, concerns then the relationship 
between these notions. We will turn to it in section~\ref{Isomorphism}.} 
Beyond them, there exists a fourth, less well-known notion of 
completeness that is discussed by Carnap in (the unpublished 
fragments of) the projected second part of his manuscript (as 
well as in his subsequent papers \citealt{Carnap1930} and 
\citealt{CarnapBachmann1936}).\footnote{What exists of part two of
\emph{Untersuchungen} is documented in Carnap's Nachlass (RC
081--01--01 to 081--01--33). See \cite{Schiemer2012,Schiemer2013}
for a detailed discussion of its contents.} That notion is called
the ``completeness of the models'' of a given axiom system by
Carnap. It is closely related to his discussion of ``extremal
axioms'', i.e., axioms (such as Peano's induction axiom in
arithmetic and Hilbert's axiom of completeness in geometry)
that impose a minimality or maximality constraint on the models of a
theory.

Extremal axioms are characterized in the following way by Carnap: 
\begin{definition}[Extremal axioms]
Let $\mP$ be a model of axiom system $f$. Then
\begin{align*}
Max(f, \mP) & \mathrel{=_{df}} \lnot \exists \mQ (\mP \subset 
\mQ \land \mP \ne \mQ \land f(\mQ)) \\
Min(f,\mP) & \mathrel{=_{df}} \lnot \exists \mQ (\mQ \subset 
\mP \land \mP \ne \mQ \land f(\mQ))
\end{align*}
\end{definition}
Here the first sentence states that $\mP$ is a maximal model of
theory~$f$: there are no elements in the model class of~$f$ that are
proper extensions of~$\mP$. Dually, the second sentence expresses
that $\mP$~is a minimal model of~$f$: there are no elements in the
model class of~$f$ that are proper submodels of~$\mP$. The
``metatheoretic'' character of Carnap's extremal axioms is again clear
enough: They do not speak about specific models of an axiom system, 
but express properties of its whole model class.\footnote{See 
\cite{Schiemer2013} for a further discussion of Carnap's extremal 
axioms, also \cite{Loeb2014} for an analysis of Carnap's notion of 
submodel underlying his work on extremal axioms.}

The type of completeness for axiomatic theories effected by a
maximality constraint---let us call it ``Hilbert completeness''---can
then be specified by Carnap as follows:
\begin{definition}[Hilbert completeness]
An axiom system $f$ is Hilbert complete iff
\[
\forall \mP \forall \mQ [(f\mP  \ \& \ f\mQ \land \mP \subseteq 
\mQ ) \rightarrow \mP = \mQ]
\]
\end{definition}
Put informally, this says that the models of theory~$f$ are all maximal, 
i.e., non-extendable to other models of $f$.\footnote{In Carnap's own 
words: ``The models of an axiom system that is closed by a maximal axiom 
possess a certain completeness property in that they cannot be extended without
violating the original axiom system.'' \citep[p.~82]{CarnapBachmann1936}. 
Carnap's main goal in part two of \emph{Untersuchungen} is to give 
an explication of different types of extremal axioms and of corresponding 
notions of completeness. As pointed out in \cite{Schiemer2013}, 
the notes for it also contain results about the relationship between 
extremal axioms and the monomorphicity, or categoricity, of a theory.}

\section{Metatheoretic notions in \emph{Untersuchungen}}
\label{Completeness}

Before turning to a more detailed assessment of Carnap's early
metatheory, several general observations concerning his definitions of
these notions are worth making. As already pointed out above, the main
difference to a current model-theoretic treatment of them concerns his
type-theoretic framework: Carnap's definitions are not formulated in a
separate metalanguage, but in a single higher-order language. An
important part of this is that the ``metatheoretic'' generalization
over models in them is expressed in terms of higher-order quantifiers
in $\Lang_{TT}$.\footnote{The question of how this squares with the
technique of model variation is discussed in section~\ref{Domain}.}

To illustrate this aspect further, consider the current model-theoretic
treatment of logical consequence. Let $\varphi$ be a sentence of a
given language $\Lang$ and let $\Gamma$ be a theory, i.e., a set of
$\Lang$-sentences. To say that $\varphi$ is a consequence of
$\Gamma$ is nowadays expressed as $\Gamma \models \varphi$ or,
more explicitly, $\forall \mathfrak{M}(\mathfrak{M} \models \Gamma
\Rightarrow \mathfrak{M} \models \varphi)$. This statement is usually
not expressed in $\Lang$ itself, but in a separate metalanguage,
e.g. the language of set theory, in which expressions of the object
language can be coded and suitable definitions of ``satisfaction'' or
``truth in a model'' given. Such an approach presupposes a strong
background theory, e.g. Zermelo-Fraenkel set theory (ZFC), in which
models of the theory~$\Gamma$ can actually be constructed. By contrast,
Carnap formulates logical consequence by means of a quantified 
conditional statement $\forall X (\Gamma(X) \rightarrow \varphi(X))$ 
in a type-theoretic language~$\Lang_{TT}$, as we saw. This statement 
does not contain model quantifiers in the now usual sense, i.e.,
quantifiers that range over the set-theoretic models of the theory in
question. Instead Carnap uses higher-order quantifiers in
$\Lang_{TT}$ that range over different possible interpretations 
of the primitive terms of the theory in question, in his sense of
interpretation.

For all its differences, there is a striking structural similarity between 
Carnap's approach and the now usual model-theoretic approach, 
at least in the case where $\Gamma$ is finite.\footnote{Note here 
that for practical purposes the restriction to finite theories is largely
irrelevant to the logicians of the 1920s. Theories we now think of
as requiring infinite axiomatizations, such as Peano Arithmetic and
ZFC, would have been given finite axiomatizations then. The induction 
schema, for instance, would be formalized in terms of a single
sentence with a higher-order quantifier. In Hilbert's axiomatics too,
the schema would be formalized as a single axiom with a formula
variable.} To see this better, one has to look closely at how
metatheoretic claims are formalized more fully today. It is a simple
model-theoretic fact that the statement ``$\varphi$ follows from
$\Gamma$'' (i.e. $\Gamma \models \varphi$) is equivalent to
``$\Gamma^{*} \rightarrow \varphi$ is valid'' (i.e. $\models
\Gamma^{*} \rightarrow \varphi$), where $\Gamma^{*}$ is the
conjunction of the sentences in $\Gamma$. The standard way to 
make the latter claim fully precise is to translate it into a set-theoretic
statement that turns out to be provable in axiomatic set theory, say in
ZFC. Thus, we can say that $\varphi$ follows from $\Gamma$ iff $\forall
M (\mathrm{Sat}(\ulcorner \Gamma^{*} \rightarrow \varphi \urcorner,
M)$ is a theorem of ZFC, i.e., iff
\[
\mathrm{ZFC} \vdash \forall M (\mathrm{Sat}(\ulcorner \Gamma^{*} 
\rightarrow \varphi \urcorner, M).
\]
Here the expression $\ulcorner \Gamma^{*} \rightarrow \varphi \urcorner$
is a coded version of the corresponding object-language statement,
$M$ is a set-variable that ranges over $\Lang$-structures, and
$\mathrm{Sat}(x,y)$ is a standard satisfaction (or truth) predicate that relates
$\Lang$-formulas to $\Lang$-structures.

Carnap's account of consequence reveals itself to be very similar to 
this set-theoretical formalization if one makes more explicit some 
assumptions concerning the \emph{Grunddisziplin} underlying his study 
of general axiomatics. Specifically, we saw in section~\ref{Axiomatic} 
that a central semantic notion taken as primitive (or left implicit) in
Carnap's account is the notion of ``being true'' or ``holding'' in
simple type theory (henceforth: TT). The truth of a theory
in a particular model is then specified in terms of this notion: a
theory (expressed by a propositional function) $f(X_1, \dots, X_n)$ is
``satisfied'' by a sequence of predicates $\langle R_1, \dots, R_n
\rangle$ if the sentence $f(R_1, \dots, R_n)$ ``holds'' in~TT. Now,
what does Carnap mean by saying that a sentence of his background
language ``holds''?

Two reconstructions would seem to be consistent with his remarks 
on this topic in \emph{Untersuchungen}. The first is to treat ``holds
in~TT'' semantically, as truth in the intended (and fixed)
interpretation of the type-theoretical language, i.e., as truth in the 
universe of types. (We will return to this interpretation in the next section.) 
The second reconstruction is to treat the notion in terms of provability in
the type-theoretic system described by Carnap. Thus, to say that $f\mR$ 
holds in the basic system can be understood as saying that $f\mR$ 
is derivable in~TT.\footnote{It is not specified in \emph{Untersuchungen} 
which basic laws and inference rules are included in TT. But Carnap 
holds at one point that ``the basic discipline has to contain theorems 
(\emph{Lehrs\"atze}) about logical, set-theoretical, and arithmetical 
concepts'' \cite[p.~61]{Carnap2000}, RC 080--34--03, p.~6.} 

Along the same lines, we can make more precise Carnap's notion of
consequence: the metatheoretic claim that $g$ is a consequence of
a theory~$f$ can be represented formally as follows:
\[
\text{TT} \vdash \forall \mR (f\mR \rightarrow g\mR)
\]
Thus, to say that $g$ is a consequence of~$f$ is simply to say that
the quantified conditional statement is a theorem of~TT. This
type-theoretic reconstruction of the consequence relation is evidently
quite similar to its now usual set-theoretic treatment. Notice, in
particular, that (first-order) ZFC and Carnap's simple type theory TT
play similar foundational roles: Both are strong background theories
that allow one to construct models of given mathematical theories. 
Moreover, the respective languages allow one to generalize
over these models in term of (first-order or higher-order) quantifiers. 
Viewed from this perspective, the main difference between 
Carnap's approach and the later model-theoretic approach is not a 
missing meta-language/object-language distinction. Rather, the 
main difference lies in the fact that the basic semantic notion of
``satisfaction in a model'' (corresponding to $\mathrm{Sat}(x,y)$ above) 
is taken as primitive, rather than recursively defined.  A recursive
definition, as model theory uses it today, would require a suitable method
of coding expressions in the meta-language, and this was unavailable to
Carnap in 1928. G\"odel and Tarski developed precise treatments of 
these notions in the following years, of course, which Carnap could
then make use of in his later work, as he in fact did.

\section{Models and domain variation}
\label{Domain}

How does Carnap's type-theoretic explication of model-theoretic notions 
square with the now usual model-theoretic approach in other respects? 
In particular, an important interpretive issue discussed in the secondary
literature concerns how the notion of domain variation is captured, and if
it can be captured at all within his framework.

Domain or model variability is a core idea in standard model theory. 
Again, a theory is expressed in a formal or disinterpreted language 
that usually contains a set of non-logical constants, the theory's
primitives. The semantic interpretation of this language is then
specified relative to a model in the usual way: The constants
are treated in terms of an interpretation function that assigns
individuals from the model's domain to individual constants, $n$-ary
relations (on the domain) to $n$-ary predicates, and $n$-ary functions
(on the domain) to $n$-ary function symbols. Quantified variables,
finally, are interpreted as ranging over the domain of the model. A
guiding idea underlying this approach is that this semantic treatment
of a language (and, consequently, of a theory expressed in it) can be
varied, i.e., it can be specified relative to models with different
domains and corresponding interpretation functions.
 
Carnap's conception of logic in the late 1920s differs clearly from
such a model-theoretic account. We already saw that theories are
expressed by him in a type-theoretic language $\Lang_{TT}$ that is 
not formal in the sense just mentioned but ``contentual''
(``inhaltlich'').\footnote{As Carnap writes: ``Every treatment
and examination of an axiom system thus presupposes a logic,
specifically a contentual logic, i.e., a system of sentences which
are not mere arrangements of symbols but which have a specific
meaning. [Jede Behandlung und Pr\"ufung eines Axiomensystems setzt
also eine Logik voraus, und zwar eine inhaltliche Logik, d.h. ein
System von S\"atzen, die nicht blo\ss e Zeichenzusammenstellungen
sind, sondern eine bestimmte Bedeutung haben.]''
\citep[p.~60]{Carnap2000}; RC 080--34--03, p.~4. As 
noted by others, \emph{Untersuchungen} does not contain 
an explanation of how the \emph{Grunddisziplin} acquires 
its specific interpretation. Such an explanation is given, at least in
outline, in Carnap's sketch \emph{Neue Grundlegung der Logik}
written in 1929. See \cite{AwodeyCarus2007} for further details.} In
other words, he uses a fully interpreted language, one whose constants
have a fixed interpretation and whose variables range over a fixed
universe of objects.\footnote{Carnap remains neutral in
\emph{Untersuchungen} with respect to the specific choice of the
signature of his background language. He holds that arithmetical and
set-theoretical terms can either be understood as ``logical''
primitives or introduced by explicit definition in the (pure)
type-theoretic language in a logicist fashion. See
 \cite[pp.~60--63]{Carnap2000}; RC 080--34--03, pp.~4--8.} Carnap is
not explicit about the precise nature of the semantics underlying
$\Lang_{TT}$ in \emph{Untersuchungen} or related writings, but---at
least in 1928---it is most likely meant to involve a rich ontology,
namely a full ``universe of types''~V. In light of his specification
of the theory of relations and of type theory in \citep{Carnap2000} 
as well as \citep{Carnap1929}, this universe of types can be
reconstructed in current terminology as a model of the form $V =
\langle \{ D_{\tau} \} _{\tau}, I \rangle$, where $\{ D_{\tau} \}
_{\tau}$ is a ``frame,'' i.e. a set of type domains, and $I$ an
interpretation function for the constants
in~$\Lang_{TT}$.\footnote{See \cite{SchiemerReck2013} for further
details, as well as \cite{Andrews2002} for a general discussion of 
the semantics of type theoretic languages.}

With this in mind, we can go back to the ``semantic'' reconstruction
of Carnap's metatheoretical approach that was mentioned as an option
in section~\ref{Completeness}. Consider again his definition of the
notion of consequence. The correctness of the statement ``$g$ follows
from theory~$f$'' can now also be represented in terms of the notion
of truth in the universe of types, or more formally, as:
\[
V \models  \forall \mR (f(\mR) \rightarrow g(\mR))
\]
According to this second reconstruction of Carnap's account, a
metatheoretic claim is correct (or holds in the \emph{Grunddisziplin}) 
if the type-theoretic sentence expressing it is \emph{true} 
in the underlying universe of types.\footnote{It is important to
emphasize that this reconstruction is most likely not how Carnap
understood the phrase ``statement $\varphi$ holds in the basic
system'' in 1928. As mentioned already, the current semantic 
notion of satisfaction was not part of his conceptual toolbox at the
time. However, and as shown by \cite{AwodeyCarus2007}, this
situation changed in the early 1930s, specifically after Carnap's
exchange with G\"odel on the notion of analyticity. In particular,
Carnap's understanding of ``metatheoretic'' statements at that point 
becomes similar to the semantic reconstruction given here, as 
argued in \cite{Schiemer2013}.}

Against that background, the question arises of whether or not, and if so
how, model domains can be varied for a given theory~$f$. To see the 
problem more clearly, recall that both~$f$ and the statement
$\forall \mR (f(\mR) \rightarrow g(\mR))$ are formulas of the language
$\Lang_{TT}$. Let us now replace $f$ in the sentence above by the
propositional function expressing a particular mathematical theory,
for instance the complex formula $\Phi_{Hausd}$ for Hausdorff's
topological theory. The resulting metatheoretic statement will then
contain two types of quantifiers that would be expressed in separate
languages today: a `metatheoretic' quantifier `$\forall \mR$',
ranging over the models of Hausdorff's theory, and `object-theoretic' 
quantifiers, used in the formulation of the axioms of $\Phi_{Hausd}$. 
(The latter are the individual and set quantifiers in Axioms 1--4 from 
section~\ref{Axiomatic}.)

For Carnap, both quantifiers are part of the vocabulary of
$\Lang_{TT}$. It follows that both come with a fixed range.  
In the case of the metatheoretic quantifiers this is intended:
Intuitively, we want such quantifiers to range over all models 
of a given theory.  Carnap's account corresponds to today's
model-theoretic approach in this respect, where model quantifiers 
are usually expressed in an informal but interpreted language 
of set theory.  The situation is different for the object-theoretic
quantifiers: As mentioned before, we usually want these
quantifiers---which are expressed in a formal object language---to 
be freely re-interpretable relative to different interpretations of
that language. Looking at Carnap's approach, it is not obvious 
in which way such model-theoretic re-interpretability of the
``object-theoretic'' quantifiers contained in $\Phi_{Hausd}$ can be
effected. Put differently, how can the model-theoretic idea of domain
variation (as assumed in concepts such as consequence 
and categoricity) be simulated in Carnap's type-theoretic framework, 
if this can be done at all?\footnote{This interpretive issue has 
been much debated in the secondary literature on Carnap's early
semantics. A central objection against his general axiomatics
project, raised early on by Hintikka, was that, due to Carnap's semantic
universalism or ``one domain assumption'', the idea of domain
variation was simply inconceivable \citep{Hintikka1991}. This view
has recently been corrected in work by Schiemer, Reck, and Loeb. In
that work, it is shown that Carnap was well aware of the importance
of capturing the notion of domain variability for his metatheoretic
work and was not without resources to do so. See
\cite{Loeb2014,Loeb2014b}, \cite{Schiemer2012,Schiemer2013}, and
\cite{SchiemerReck2013}.}

One way to simulate this kind of domain variation that was frequently 
employed by Carnap at the time consists in the method of quantifier
relativization. The idea is to effectively restrict the range of an
object-theoretic quantifier used in the formulation of a theory to the
interpretation of the primitive terms of the theory. Thus, while
object-language quantifiers have a fixed interpretation in Carnap's
account in general, their range can be relativized to a particular
``model domain'' of a theory as soon as that theory is ``interpreted''
in a model in Carnap's sense. In fact, one can distinguish between two
ways of relativizing quantifiers to model domains of theories along
such lines: (i) the direct relativization via specific ``variabilized'' 
domain predicates (which itself comes in two sub-variants) and (ii) 
the indirect relativization via other primitive terms of a theory. It is 
instructive to look at some concrete examples from Carnap's writings 
so as to understand better how these methods are used by him.

In some cases, Carnap restricts the quantifiers used in 
the formulation of an axiomatic theory in terms of a specific 
unary ``domain predicate'' that is introduced with the primitive 
terms of that theory (either as itself a primitive or as a defined
predicate).\footnote{As pointed out by \cite{Schiemer2013}, a very
similar convention of type relativization can also be found in Tarski's 
work on the ``methodology of the deductive sciences,'' from the 
same period. Compare \cite{Mancosu2010} for a detailed discussion 
of Tarski's case and the other secondary literature on the topic.} A 
typical example is the formalization of Peano arithmetic in \emph{Abriss}. 
The ``original'' axiomatization presented there is based on three primitive 
terms, namely $\mathsf{nu}$, $\mathsf{za}$, $\mathsf{Nf}$ (standing for
``Zero'', ``$x$ is a number'' and ``the successor of $x$'', respectively). 
By looking at Carnap's formulation of the corresponding axioms, 
it becomes clear that in this case the predicate $x \in \mathsf{za}$ 
functions as a domain predicate in the way described above:
\begin{align}
\tag{PA1} & \textsf{nu} \in \mathsf{za} 
\\
\tag{PA2} \forall x & ( x \in \mathsf{za} \rightarrow \mathsf{Nf}(x) \in \mathsf{za}) 
\\
\tag{PA3} \forall x\forall y &( ( x,y \in \mathsf{za} \land \mathsf{Nf}(x) = 
\mathsf{Nf}(y)) \rightarrow x = y )
\\
\tag{PA4} \forall x &( x \in \mathsf{za} \rightarrow \mathsf{Nf}(x) \neq \mathsf{nu})
\\
\tag{PA5} \forall \alpha & (( \mathsf{nu} \in \alpha \land \forall x ( x \in \alpha 
\rightarrow \mathsf{Nf}(x) \in \alpha)) \rightarrow  \mathsf{za} \in \alpha)
\end{align}
In other examples from Carnap's work on axiomatics such domain
predicates are not part of the theory's signature, but are introduced
by means of explicit definitions from the primitive vocabulary.\footnote{This 
approach of specifying domain predicates was first discussed by  
\cite{Loeb2014}.} Consider again Carnap's formulation of Hausdorff's
topological axioms introduced in Section~\ref{Axiomatic}: Here a set
variable $\mathsf{pu}$ (standing for ``$x$ is a point'') is introduced
by definition, based on the primitive term $\mathsf{U}(x,y)$ (standing
for `$x$ is a neighborhood of $y$'). The individual quantifiers needed
in the formulation of the axioms are then, in the relevant cases, 
relativized to~$\mathsf{pu}$. Take axiom (Ax4) as an illustration:
\begin{align}
\tag{H4} & \forall x, y (x, y \in \mathsf{pu} \land x  \neq y \rightarrow 
\exists \alpha, \beta ( \mathsf{U}(\alpha,x) \land \mathsf{U}(\beta,y)  
\land \alpha \cap \beta = \emptyset))
\end{align}
Finally, there are several examples of axiom systems in Carnap's work
where unary domain predicates (or set predicates) are not employed
at all. Instead, the quantifiers used in the formulation of a theory are
relativized to model domains more indirectly and implicitly, in terms
of other primitive terms of a theory. A typical example is a second
axiomatization of ``basic arithmetic'' (BA) discussed in
\emph{Untersuchungen} (also, later, in
\citealt{CarnapBachmann1936}). Here the axioms of BA involve a single
binary predicate $R(x,y)$ (for ``$x$ is the successor of $y$''), as
follows:\footnote{Carnap and Bachmann add a fourth meta-axiom to
  BA1--3 that restricts the possible interpretations of BA1--3 to
  minimal models, in the sense specified in
  Section~\ref{Completeness}.  See
  \citep[p.~179]{CarnapBachmann1936}.}
\begin{align}
\tag{BA1}  \forall x \forall y &(R(x,y) \rightarrow \exists z(R(y,z)) )
\\
\tag{BA2} \forall x \forall y \forall z &( (R(x,y) \land R(x,z) \rightarrow y=z) 
\land (R(x,y) \land R(z,y) \rightarrow x=z))
\\
\tag{BA3} \exists! x & (x \in \mathrm{Dom}(R) \land x \notin \mathrm{Ran}(R))
\end{align}
Notice that in this case the object-theoretic quantifiers are effectively 
restricted to the field of any binary relation assigned to the primitive 
sign $R(x,y)$, and thus, to any given model domain of the theory~BA.

Given these examples, two additional remarks about Carnap's 
approach are in order. First, one main difference between his
approach and current model theory is that instead of providing an
interpretation function for un-interpreted non-logical symbols,
Carnap's approach simply quantifies over these non-logical
symbols. The primitive terms of a theory are expressed by free
variables; and the semantic specification of the latter, relative to 
a particular model, is given in terms of the substitution of these
variables by interpreted constants of the basic discipline.  It might
seem as if this way of doing things leaves no room for another crucial
ingredient of model theory, namely the domain of quantification
of a model.  But as we have seen, this is accommodated in Carnap's
account: domains are understood either as the extensions of a specific
(primitive or defined) domain predicate or, in cases where such 
a predicate is absent, as the (first-order) fields of the relations
assigned to the primitive signs of a given theory.\footnote{See
\cite{Schiemer2013} for a more detailed discussion of Carnap's
``domain as fields'' conception of models. Compare \cite{Loeb2014}
for an alternative account of Carnap's understanding of
models. Loeb's discussion focuses on examples of axiomatic theories
where a domain predicate is introduced into the language in terms of
an explicit definition from the theory's primitive signs. These
defined domain predicates provide further confirmation for the
interpretation of Carnap's conception of model given in
\cite{Schiemer2013}, given that in most cases a model domain is
explicitly specified as the domain or range or field of a given
primitive relation. Nevertheless, Loeb holds that ``the
domains-as-fields conception is too strict to describe Carnap's
practice.''  \cite[p.~427]{Loeb2014} The particular example she has
in mind is an axiomatization of projective geometry with one unary
primitive predicate $ger$ (for the class of `lines'). The domain of
a model of the theory (i.e. the class of `points') is not defined as
the field of the relations assigned to \emph{ger}, but as the union
of the elements of all lines \cite[\S 34]{Carnap1929}. We would like
to stress that this and similar examples are fully in accord with
the interpretation of model domains given in
\cite{Schiemer2013}. For a more detailed analysis of model domains
of lower types see, in particular, the discussion of Carnap's notes
on ``domain analysis'' from part two of \emph{Untersuchungen} in
\cite[pp.~506-508]{Schiemer2013}.}

Second, it should be clear by now that the various ways of quantifier
relativization used by Carnap allow him to simulate the domain
variability typically assumed in model-theoretic notions (such as
consequence and categoricity) despite the fact that his theories are
expressed in a single fully interpreted language. Instead of
reinterpreting a formal object language (including the range of the
object-language quantifiers) in a separate metatheory, domain
variation is effected by the relativization of the quantifiers of his
single language to the primitive terms of a theory. Carnap's
systematic use of this technique shows (contra Coffa, Hintikka,
and others) that his explications of various notions were genuinely
model-theoretic in spirit. His proposal contained an early form of
model theory that included the notion of domain variation.

\section{Consequence and derivability}
\label{Consequence}

In modern metalogic, we use two important and related
distinctions. The first, already discussed above, is that between the
object language, for which we give definitions of, say, logical
consequence and derivability, and the metalanguage, in which we
formulate these definitions and prove corresponding results.  The
second distinction, to be discussed further now, is that between
semantic notions, such as truth, consequence, and satisfiability, and
proof-theoretic notions, such as proof, derivability, and consistency.
One of the main tasks of metalogic is to relate these different notions, 
e.g., when proving soundness and completeness. Even before Tarski 
made the object/metalanguage distinction a topic of explicit study, 
logicians respected that distinction implicitly (at least to a large extent).
For example, Carnap's \emph{Untersuchungen} is written in German
(the meta-language), in which definitions for his type-theoretic
language (the object-language) are formulated and corresponding 
results proved. As will be clarified further in the next section, Carnap 
was also aware of the distinguishing between truth and proof. Yet 
both distinctions were neither as sharply drawn nor their importance 
as fully recognized by him as one would expect today.

As just noted, Carnap formulated his relevant definitions in a
metalanguage. But instead of defining notions such as consequence 
and satisfiability by means of metalinguistic notions, he 
defined them by using sentences of his \emph{Grunddisziplin}, 
thereby taking for granted that its sentences express determinate 
propositions.  For instance, recall the definition of satisfiability of 
an axiom system~$f$, which we encountered in passing above: 
``$f$ is satisfiable iff $\exists \mR\, f\mR$.'' From today's viewpoint, 
this definition fails to distinguish between expressions of an object 
language and expressions that now belong to a separate metalanguage, 
such as the quantifier ``$\exists \mR$'' and predication of the relation 
``satisfies'' to $\mR$~and~$f$.  In addition, Carnap talks about ``values''
of the variables~$\mR$---his models, which can be exhibited and 
investigated in the metalanguage. It is in the latter connection that 
his approach led to several further confusions and false starts, as a 
closer look at \textit{Untersuchungen} will reveal.

One problem Carnap did not appreciate in this context is that
describing a model of $f$ in his metalanguage and proving $\exists
\mR\, f\mR$ in type theory are by no means the same thing.  That is to
say, a metatheoretic description of a model of $f$ does not, by
itself, result in terms in type theory that can be substituted for
$\mR$ in $f\mR$ so that $f\mR$, and with it $\exists \mR\,f\mR$, can
be proved from the axioms of type theory. Conversely, type theory may
prove $\exists \mR\,f\mR$ without a description of the model being
available in the metalanguage. Having said that, Carnap's definition
of satisfiability is acceptable if two additional assumptions can be
made: that all relevant mathematical constructions which can be
carried out in the metalanguage can also be formalized in type 
theory; and that what type theory proves is true. It is not implausible 
to suppose that Carnap took these assumptions for granted.  The first
assumption received ample evidence from the success of Russell and
Whitehead's work in formalizing mathematics, the second from the
belief that the axioms of type theory express logical truths (at least
those of the simple theory of types).

A second problematic aspect, already mentioned briefly above, is that
the satisfaction relation is not explicitly defined by Carnap. Indeed,
instead of taking it to be a relation in the metalanguage between the
axiom system $f$, seen as a syntactic object, and the model $\mR_1$,
seen as a semantic object, and then using a recursive definition along
Tarskian lines, he simply assumes that ``$\mR_1$ satisfies $f$'' is
meaningful. In section 2.3 of
\emph{Untersuchungen},\footnote{\citet{Carnap2000}, p.~95; RC
  080--34--02, \S13, p.~42.} Carnap defines an ``admissible model'' of
$f\mR$ as a system of relations $\mR_1$ consisting of ``allowed
values'' of the variables in~$\mR$, i.e., a system of relations of the
same type as the types of the variables in~$\mR$. If it is, he writes,
then the propositional function $f\mR$ applied to $\mR_1$ does or
doesn't ``hold''. This is imprecise by today's standards, although it
can be interpreted in a way that it is correct, as we saw. If $\mR_1$
is an admissible model of $f\mR$, Carnap also sometimes writes that
$f\mR_1$ is ``true'' or ``false''. This is worse, since ``$\mR_1$'' is
not an expression of type theory, and hence, neither is $f\mR_1$.

To some degree we may forgive Carnap this lapse, since everyone at 
the time spoke this imprecisely. For instance, logicians in the Hilbert 
school also considered basic semantic notions, such as satisfaction 
and truth, and spoke freely of formulas of first order logic being
``satisfied'' in certain domains, or of a formula being ``true'' for a 
particular value of its free variables, as if sets and relations could all
automatically serve as names for themselves in the object language.
Moreover, when satisfaction in a particular model is not at issue but
only general claims about satisfaction and models, this problem does
not arise: all the model variables are then bound by type-theoretic
quantifiers and the statements in question are sentences of type
theory. In that sense (and assuming that the sentences of type theory
are meaningful and express logical propositions), Carnap was justified
in regarding his results to be theorems.

A third crucial difficulty in Carnap's definitions of metalogical concepts
concerns his quantification over propositional functions.  For
instance, Carnap's definition of ``$f$ is inconsistent'' is given by
the following type-theoretic sentence:
\[
\exists h\forall \mR(f\mR \to (h\mR \land \lnot h\mR))
\]
When defining ``$f$ is satisfiable'' in terms of $\exists \mR\, f\mR$,
this at least expresses what is intended (the existence of a
satisfying sequence of relations of the right types). By contrast, 
the definition here does not.  For $h$ is a variable ranging 
over higher-type objects, it seems, and not, as would be needed,
a variable ranging over syntactic expressions.  Again, Carnap is not
alone in being unclear about this aspect at the time.  Just as in
\emph{Principia Mathematica}, his $h$ is a variable ranging over
propositional functions; and Russell and Whitehead were notoriously
vague about whether quantification over propositional functions was
supposed to be thought of substitutionally or objectually.  But there
are more profound problems with Carnap's approach to consistency.

Today we use the terms ``consistent'' and ``inconsistent'' as
proof-theoretic notions: an axiom system is inconsistent if there is a
sentence $\varphi$ such that $\varphi \land \lnot \varphi$ can be proved 
from it. Carnap's definition, by contrast, uses his own notion of
consequence: an axiom system is inconsistent if there is a sentence
$\varphi$ such that $\varphi \land \lnot \varphi$ is a consequence of it.  Now,
this was an idiosyncratic definition even when Carnap was writing.
We know that Hilbert, who put great emphasis on the consistency of
axiom systems (and proofs thereof), already defined it proof-theoretically
at that time.\footnote{However, the usage of logical terminology was by no
means settled then. For example, \cite{Skolem1920} used the term 
``inconsistent [widerspruchsvoll]'' to mean unsatisfiable; even Hilbert himself 
at times uses ``$A$ is a consequence of [\emph{folgt aus}] $B$'' to mean 
``$A \to B$ holds.''}  Carnap was certainly aware of Hilbert's work and,
one would assume, intended to provide a definition that matched Hilbert's.
The fact that he thought his definition managed to do that seems to be
the result of several further confusions.  

To begin with, Carnap thought that provability in Hilbert's sense and
consequence as defined by him coincide.  In fact, he presented an
argument for this: Section 2.2 of \textit{Untersuchungen} contains an
explicit comparison of his approach with
Hilbert's.\footnote{\cite{Carnap2000}, p.~90ff; RC 080--34--03, \S12,
p.~39ff.}  Carnap establishes there that, if $gx$ can be proved from
the axiom system $fx$ in Hilbert's sense, then $\forall x(fx \to gx)$
can be proved, and conversely.  He concluded that the concept of
consequence defined by ``is provable from'' and his own notion of
consequence coincide.  But this inference is unwarranted, for two
reasons.  One is technical: Carnap's notion of consequence relies on a
much stronger system than Hilbert's, and, e.g., the provability of
$\forall \mR(f\mR \to g\mR)$ in type theory does not obviously
guarantee the provability of the formula $f\mR \to g\mR$ in Hilbert's
system.\footnote{If the variables in $\mR$ are at most second order,
then the result holds.  In that case, assuming soundness of type
theory, $f\mR \to g\mR$ is valid in first-order logic. Thus its
provability in Hilbert's system follows from G\"odel's completeness
theorem, also from a proof-theoretic conservativity proof, neither 
of which were available to Carnap yet.}  The other reason is more
fundamental: Carnap's argument establishes that $g\mR$ is provable
from $f\mR$ in Hilbert's system iff $\forall \mR(f\mR \to g\mR)$ is
provable in type theory.  But that is not enough to establish that
Carnap's consequence coincides with Hilbert's provability.  For that,
we would need to show that $g\mR$ is provable from $f\mR$ iff 
$\forall \mR(f\mR \to g\mR)$ is \emph{true} in type theory.

In section 2.4 of \emph{Untersuchungen}, Carnap considers the
relationship between consistency and satisfaction.  In particular, he
argues that every inconsistent axiom system is unsatisfiable (Theorem
2.4.1), and conversely, that every consistent axiom system is
satisfiable (Theorem 2.4.8).\footnote{\S2.4 of
\cite[p.~96ff]{Carnap2000} is \S14 in RC 080--34--03, p.~46ff;
Theorem 2.4.1 is Satz~1 on p.~47; Theorem 2.4.8 is Satz~8 on p.~49.}
The definitions at work here are Carnap's, of course, not the ones
used today.  Thus, Theorem 2.4.1 establishes, in type theory, that
\begin{align}
\exists h\forall \mR(f\mR \to (h\mR \land \lnot h\mR)) & \to \lnot\exists \mR\, f\mR, \label{snd}\\
\intertext{and Theorem 2.4.2 the reverse direction, in the form of its contrapositive,}
\lnot\exists h\forall \mR(f\mR \to (h\mR \land \lnot h\mR)) & \to \exists \mR\, f\mR. \label{cmpl}
\end{align}
Carnap's proofs consist of several elementary logical steps.  For
instance, to establish (\ref{cmpl}), he argues as follows: Assume
the antecedent, $\lnot\exists h\forall \mR(f\mR \to (h\mR \land \lnot
h\mR))$.  Pushing the initial negation across the quantifiers and
using $\lnot(\varphi \to \psi) \to (\varphi \land \lnot\psi)$, we obtain
$\forall h\exists \mR(f\mR \land \lnot (h\mR \land \lnot h\mR))$.
As $f$ only appears in the left conjunct and $h$ only in the right,
we can shift quantifiers and obtain $\exists \mR\,f\mR \land \forall
h\lnot (h\mR \land \lnot h\mR)$.  The left conjunct is the desired
conclusion.  So indeed Carnap's proofs are correct---but the results 
do not mean what we mean by ``consistency implies satisfiability''.

Carnap's argument for the equivalence of Hilbert's notion of
provability and his own notion of consequence does yield the
following: If an axiom system is inconsistent, i.e., can prove $\varphi
\land \lnot \varphi$ for some formula $\varphi$, then type theory proves the
antecedent of Theorem 2.4.1, ``$\exists h\forall \mR(f\mR \to (h\mR
\land \lnot h\mR))$,'' and hence also the consequent.  But the negation 
of ``$f$ is (Hilbert) inconsistent'' is only: for no formula $\varphi$
does $f$ prove $\varphi \land \lnot \varphi$.  It does not follow that type
theory proves $\lnot\exists h\forall \mR(f\mR \to (h\mR \land \lnot
h\mR))$ (and Carnap does not give an argument for the latter).  In
fact, by G\"odel's incompleteness theorem we know that this is false 
when the axiom system includes second-order quantifiers: take $f$ 
to be the conjunction of the (finitely many) second-order axioms of
$\mathit{PA}^2$ together with a false-in-$\mathbb{N}$-but-not-refutable 
sentence, such as $\lnot\textrm{Con}(\mathit{PA}^2))$.  This 
sentence is consistent but unsatisfiable, as G\"odel showed.  Hence 
$\exists \mR\,f\mR$ is false, and so is 
$\lnot\exists h\forall \mR(f\mR \to (h\mR \land \lnot h\mR))$, 
since the two are equivalent, as Carnap showed.  So,
while provability (in particular, inconsistency) in our sense implies
provability of consequence (and inconsistency) in Carnap's sense,
unprovability (and consistency) in our sense do not imply provable
consequence (provable consistency) in Carnap's.  In fact, full type
theory is strong enough to prove $\lnot\exists \mR\,f\mR$.

As we have seen, there are a number of unclarities, gaps, and other
problems in Carnap's approach to derivability.  At this point, the
question arises whether what he did can be interpreted charitably, or
can be reconstructed so that it becomes correct.  Compare here 
what G\"odel writes in the introduction to his dissertation:
\begin{quote}
If we replace the notion of logical consequence (that is, of being
formally provable in finitely many steps) by implication in Russell's
sense, more precisely, by formal implication, where the variables 
are the primitive notions of the axiom system in question,
then the existence of a model for a consistent axiom system (now taken
to mean one that implies no contradiction) follows from the
fact that a false proposition implies any other, hence also every
contradiction (whence the assertion follows at once by indirect
argument). \citep[p.~63]{Godel1929}
\end{quote}
What G\"odel says here is correct, of course, but it is not quite what
Carnap showed (G\"odel credits the observation to Carnap): for
$\exists h \dots$ still does not---as pointed out above---quantify
over the possible consequences of the axioms, i.e., over formulas.
And even when the quantifier is given a substitutional reading, the
possible values of $h$ are all expressions of type theory with only
the variables in $\mR$ free.  (This is relevant also in the case of
the \emph{Gabelbarkeitssatz}, see section~\ref{Isomorphism} for
more on this point.)

Overall, Carnap's attempt to characterize provability in terms of
(provability in type theory of) a corresponding sentence expressing
consequence fails.  Once more, Carnap mistakenly assumed that
unprovability (consistency) implied provability in type theory of the
negation of the corresponding sentence.  He also mistakenly assumed
that what should be a meta-level quantification over formulas can be
expressed on the object level---in type theory---by a quantification
over propositional functions.  Finally, there is a further difference between
Carnap's approach and the now common one, a difference that relates to 
the implicit meta-level quantification over axiom systems.  When proving
completeness, say, we prove it for all axiom systems; and when such 
a proof is formalized in a formal system (say, ZFC), the formalization 
contains a quantifier over representations of axiom systems (or at
least, representations of single axioms).  Carnap never does that---and 
this is noteworthy also for another reason.  He could have expressed his 
``completeness'' theorem (that consistency implies satisfiability, in his 
sense of the terms) in a single sentence of type theory, viz.,
\[
\forall f(\lnot\exists h\forall \mR(f\mR \to (h\mR \land \lnot h\mR)) 
\to \exists \mR\, f\mR)
\]
Moreover, his proof for it would have gone through.  The fact that he 
did not do that, despite the fact that he had no qualms about formalizing 
the quantification over potential contradictions~$h$ in general, suggest 
that he might have had an inkling that such a quantification is problematic.

\section{Isomorphism and Carnap's notion of ``formality''}
\label{Isomorphism}

One important detail of Carnap's \emph{Untersuchungen} that has not
yet received significant attention in the literature is Carnap's notion of 
formality (already mentioned above in passing).  Carnap motivates the 
need for this notion in Section~2.5.\footnote{\S2.5 of \cite{Carnap2000} 
is \S15 in RC 080--34--04, p.~52ff.}  In order to properly capture 
independence from an axiom system, it must be possible to restrict the
propositional functions considered (as being a consequence of, as being
inconsistent with, or as being compatible with the axiom systems).  The
reason is that propositional functions in general can make reference
to elements of the \emph{Grunddisziplin}.  For instance, a
propositional function $h\mR$ may express ``the field of $\mR$
contains the number~20.''  Then for any satisfiable axiom
system~$f\mR$, both $h\mR$ and $\lnot h\mR$ are compatible with
$f\mR$, since any model of~$f\mR$ which contains $20$ in the field of
its relations is isomorphic to one that does not.  Hence, if no
restriction is put on the propositional function~$h$, every axiom
system is ``forkable'' by such an~$h$.  In addition, an axiom system
itself may make reference to the elements of the fields of its
relations.  As a consequence, such an axiom system~$f$ may be
satisfied by some model $\mR_1$ but not by another model $\mR_2$ even
when $\mR_1$ and $\mR_2$ are isomorphic.

Carnap's main aim in \emph{Untersuchungen} was a proof of his
\emph{Gabelbarkeitssatz}: every monomorphic (categorical) axiom system
is non-forkable, and conversely.\footnote{See \cite{AwodeyCarus2001}
  and \cite{AwodeyReck2002} for further background.}  Even for just
the left-to-right direction of this equivalence to hold, the scope of
axioms systems must be restricted to those that do not make reference
to the specifics of the domain they are describing (what Carnap calls
their ``Bestand''); and similarly for potential forking sentences.  As
formalized today, this is done by restricting the language in which
the axiom system, on the one hand, and its potential consequences, on
the other hand, are formulated.  Carnap's definitions and results are
schematic, i.e., he gives instructions for how to express, e.g., that
$f\mR$ is monomorphic.  But the sentence in type theory that expresses
that fact depends for its form on~$f\mR$.  In particular, this
sentence quantifies over the variables appearing in~$f\mR$, and so the
quantifier prefix changes depending on the number, order, and types of
the variables in~$f\mR$.  Such definitions could in principle be
restricted to propositional functions~$f\mR$ that syntactically do not
contain certain elements, such as constants of the
\emph{Grunddisziplin}.  However, wherever Carnap needs quantification
in type theory over propositional functions, this is no longer
possible. For example, in the definition of ``$f$ is forkable,'' 
\[
\exists h(\exists \mR (f\mR \land h\mR) \land \exists \mS (f\mS \land
\lnot h\mS) \dots),
\]
the quantifier $\exists h$ cannot be restricted to only $h$'s of a
certain form, because it does not actually quantify over expressions
but over higher-type objects.

Carnap's response to this problem is to define a property of axiom
systems and of their potential consequences meant to capture the idea
that they should only specify ``structural'' properties of their
models, as opposed to ``contentual'' properties such as whether $20$
is a member of the domain or not.  Given a propositional function 
$f\mR$, a value of $\mR$ of the correct types is an ``admissibe model;'' 
and it is a ``model of $f\mR$'' if it is admissible and satisfies~$f\mR$.  
Carnap now defines a propositional function to be ``formal'' if every 
admissible model that is isomorphic to a model of $f\mR$ is also 
a model of~$f\mR$. In type theory, this can be expressed by
\[
\For(f) \equiv \forall \mP\forall \mQ((f\mP \land \Iso_q(\mP, \mQ)) \to f\mQ),
\] 
(see p.~124).\footnote{RC 080--34-03, p.~75.}  Here $\Iso_q(\mP,\mQ)$ is
given by Carnap's definition of $q$-level isomorphism of the models
$\mP$ and $\mQ$. (Carnap used this notion to define isomorphism of
models for arbitrary sequences of variables of arbitrary types, as
discussed further below.)

Carnap's resulting restriction of his central theorems, especially the
\emph{Gabelbarkeitssatz}, to formal axiom systems and formal
consequences is essential to the proofs he gave for them.  In fact,
they make these proofs almost trivial.  For instance, Theorem 3.4.3 in
\emph{Untersuchungen}, that every forkable axiom system is polymorphic
(i.e., not monomorphic), is proved as follows:\footnote{\S26, Satz 3,
RC 080--34--03, p.~86.}  An axiom system~$f$ is forkable if there is
a formal~$g$ such that both $g$ and $\lnot g$ are compatible with
$f$. That is to say, there is a model $\mP$ of $f \land g$ and a model
$\mQ$ of $f \land \lnot g$.  In other words, $\mP$ and $\mQ$ are both
models of $f$, $\mP$ is a model of $g$, and $\mQ$ is a model of $\lnot
g$.  Because $g$ is formal, if $\mP$ and $\mQ$ were isomorphic, since
$\mP$ is a model of $g$, $\mQ$ would also have to be a model
of~$g$. And since $\mQ$ is a model of $\lnot g$, $\mP$ and $\mQ$ 
are not isomorphic.  The crucial fact used here is that, because $g$ is 
formal, if it is true in one model it is also true in any isomorphic one.
Note that this fact needs no proof here; it simply follows from the
definition of ``formal.''

Theorem 3.4.2, the converse direction of the \emph{Gabelbarkeitssatz},
is more intricate.\footnote{\S26, Satz 2, RC 080--34--03, p.~86.}  It
states that every polymorphic axiom system is forkable.  In its standard 
formulation today---against which the ``correctness'' of Carnap's
result is often measured---this claim depends crucially on the
language in which the axiom system is formulated.  If the language is
first-order, the corresponding claim is false: because of the
L\"owenheim-Skolem Theorem, every first-order theory with infinite
models has non-isomorphic models, even when that theory is
semantically complete (i.e., when every sentence in the language is
either a consequence of or inconsistent with it).  But working along 
Carnap's lines, there is no clear way in which the language of the 
axiom system can be restricted to first-order.  While we can consider a
theory with only first-order relations, the axioms themselves may
always contain quantified variables of arbitrary types.  Once more,
Carnap's definition of an axiom system's being forkable is given in 
type theory itself:
\[
\exists h(\For(h) \land \exists \mR (f\mR \land h\mR) \land \exists
\mS (f\mS \land \lnot h\mS))
\]
However, as was indicated by \citet{AwodeyCarus2001} and
\citet{AwodeyReck2002}, the truth of this sentence does not imply the
existence of a propositional function with only the variables~$\mR$
free.  Moreover, to prove it in type theory, it is not required to exhibit
such a propositional function, but only to show in type theory that a
value of~$h$ with the requisite properties exists.

Carnap's actual proof of this questionable direction of the
\emph{Gabelbarkeitssatz} is a perfectly valid line of reasoning in
type theory.  He did, however, misunderstand and overestimate its 
significance. Carnap begins by assuming the antecedent of the 
conditional to be proved, i.e., that there are non-isomorphic 
models of $f$:
\[
\exists \mP\exists \mQ(f\mP \land f\mQ \land \lnot\Iso(\mP, \mQ))
\]
He then picks $\Iso(\mR,\mP)$ as the required~$h$.  If we insist on
reading Carnap as attempting to prove a result for first-order logic,
we would locate his mistake here: $\Iso(\mR,\mP)$ is not a forking
sentence of the required kind. That is, $\Iso(\mR,\mP)$ is not a
formula in which the only free variables are $\mR$. The problem is
with the $\mP$: Either we think of it syntactically, as a (sequence
of) variable(s).  But these occur in the scope of the quantifier
$\exists \mP$, and so do not have the same function when considered
outside this context. Or, we might take Carnap to be arguing
semantically. In that case, $\mP$ is not a sequence of bound
variables, but a value of such a sequence of variables and so
not properly ``part of''~$h$.  However, if we work in type theory, as
Carnap does, there is no need to isolate the propositional
function~$h$. If there is a model $\mP$ of $f$, then there is a class
of models isomorphic to~$\mP$, and this is an admissible value of the
higher-type variable~$h$.  This value satisfies~$\For(h)$, and that
fact can be proved in type theory just from the definition of $h$ and
$\For$, without making use of the definition of~$\Iso_q$. As $\mP$ is
isomorphic to itself, and $\mQ$ by assumption is not, $h\mP$ and
$\lnot h\mQ$; so the sentence expressing that $f$ is forkable follows.
Note that the notion of isomorphism does almost no work here. It is 
needed only to prove that $\mP$ is isomorphic to itself, which is 
used in the last step, as well as to establish that $h$ is formal---but 
since formality is defined in terms of isomorphism, that is also trivial.

We mentioned before that Carnap's restriction on the definition of
forkability relieves him of the need to actually prove that isomorphic
models satisfy the same propositional functions: this is because he
only considered formal propositional functions, for which this property 
holds by definition.  But then, the direction of the \emph{Gabelbarkeitssatz} 
that remains true when re-interpreted using current definitions---namely 
that isomorphic structures are equivalent in the sense of satisfying the 
same propositional functions---cannot directly be ``extracted'' 
from Carnap's proofs.  Indeed, this claim is false given Carnap's setup, 
because in the language of the \emph{Grunddisziplin} a propositional 
function may require its models to contain specific objects, e.g., the 
number $20$, and so won't be satisfied by a model that does not, 
even if that model is isomorphic to a model that satisfies it.  

In order to make the result true without the restriction to formal 
propositional functions (which, again, makes the result true but trivial), 
one would have to restrict the language of the \emph{Grunddisziplin}.  
This is possible and, in fact, was done by \citet{LindenbaumTarski1935}.  
In contrast to Carnap, they were very clear on the distinction between 
the correct result,
\begin{align*}
\forall \mP\forall \mQ(\Iso(\mP, \mQ) & \to h\mP \equiv h\mQ)\\
\intertext{where $h$ is considered a schematic variable for 
propositional function expressions with the right kind of restrictions, 
and the incorrect}
\forall h\forall \mP\forall \mQ(\Iso(\mP, \mQ) & \to h\mP \equiv h\mQ)
\end{align*}
As Lindenbaum and Tarski point out, the latter is not only not provable 
but refutable. Why is that so?  Because, again, type theory proves the 
existence of a value of~$h$ for which it is false: Take the class of all 
models of the right type which contain $20$ in the field of one of its 
relations.

Neither direction of Carnap's proofs of the \emph{Gabelbarkeitssatz}
actually made substantive use of the general notion of isomorphism
that Carnap defined.  Nevertheless, it is interesting to take a look at 
what his definition of isomorphism is.  Although both the term and the
basic idea behind it predate Carnap's discussion, it is instructive to 
reconsider his attempt at a general definition for two reasons: First, 
it is, as far as we know, the first attempt at defining this notion for higher 
type structures; and second, it gives further insight into the notion of
model at work in his writings.

The notion of isomorphism (and the crucial metalogical property 
based on it, namely that of categoricity) goes back at least to
\cite{Veblen1904}.\footnote{Compare \cite{AwodeyReck2002} and the
references in it.}  Russell and Whitehead defined isomorphisms between 
relations in \emph{Principia Mathematica}, and the same approach 
can be applied to the relations that constitute Carnap's models of axiom 
systems.  This is, indeed, the insight from which Carnap proceeds.
A relation isomorphism~$S$ between two relations $R$ and $R'$, in 
the Russell-Whitehead sense, is a relation which is one-one (i.e.,
bijective) and which has the property that whenever $aRb$ there are
$a'$ and $b'$ so that $aSa'$, $bSb'$, and $a'R'b'$, as well as
conversely.  Note here that, in the Russell-Whitehead theory, a relation
$R$ is given just by its graph: it is not specified as a relation on some 
domain.  Consequently, there can be no object that is part of the
``domain'' of $R$ but is not related to any object.  Moreover,
isomorphisms are not required to relate any objects that do not
participate in the relations $R$ and $R'$.

If, as Carnap did and as was common at the time, one takes relations
in this sense to be the interpretations of the non-logical constants
of an axiom system, it is never necessary to also specify an
individual domain for the models of the axiom system: a model is
simply all the relations that interpret the individual non-logical
constants taken together.  This is, of course, closely related to the
issue discussed above about domain variability in Carnap's general
axiomatics.  Indeed, for none of the axiom systems considered in
Carnap's \textit{Abriss}, for instance, would it be necessary to
entertain the possibility of objects that are not in the fields of
their relations. When quantification over the domain is necessary, it
is formalized using a relativized quantifier, where the relativization
clause uses either a primitive non-logical constant or a defined one,
e.g., restriction to the field of a primitive non-logical relation.

As Carnap wanted to provide a general theory of axiom systems in
\emph{Untersuchungen}, he needed a general definition of when two
models are isomorphic.  No such general definition had been given
before.\footnote{Typically isomorphisms of models was discussed only
for models of a particular axiom systems, and often these only
contained a single non-logical primitive. For instance,
\citet{Veblen1904} considered a version of geometry where the single
primitive is a three-place collinearity predicate relating points.
\cite{vonNeumann1925} also discussed the (non)categoricity of his
axiom system; but his axioms were formulated in German and both what
counts as a model and what counts as an isomorphism was left vague.}
Hence Carnap had to provide one himself.  For axiom systems involving
only one non-logical relation, the notion of relation isomorphism mentioned
above suffices.  But in the absence of a domain specification as part
of the model it is not straightforward to define a notion of model
isomorphism for multiple relations.  If these relations are all just
first-order, it suffices to require the isomorphism~$S$ to preserve
all the relations.  Carnap does not do this directly.  Instead, he
only requires that the relations of the two models are pairwise
isomorphic, as well as that the individual relation isomorphisms 
agree on the intersection of the fields of those relations.

Carnap may have intended this only as a first stab at a definition, 
to be improved later on; but it also reflects a limitation of the type
system with which Carnap is working (which, admittedly, does not yet
appear if all relations are first order).  For instance, if the axiom
system contains a one-place first-order predicate constant $P_{(0)}$
and a one-place second-order predicate constant $Q_{((0))}$, then the
values of the corresponding variables are a class of individuals and a
class of classes of individuals.  But because types are not cumulative
in Carnap's system, no relation can be a bijection with a field that
contains both individuals and classes, i.e., no single $S$ can be a
relation isomorphism between both $P$ and $Q$.  

A natural solution to this difficulty caused by Carnap's type system is to 
require one isomorphism for $P$ and another for $Q$.  Of course, this
will not do yet, as Carnap explains himself.  Merely relating the classes 
falling under $Q$ in one model with those falling under $Q$ in the other 
model is not enough: the classes themselves must be isomorphic.  Carnap 
accomplishes this in \emph{Untersuchungen} by using his definition of a 
$q$-correlator, i.e., a relation isomorphism for relations of $q$-th order.  
Moreover, the isomorphisms for the individual variables must also match 
up with each other.  Carnap accomplishes the latter by requiring a 
$q$-isomorphism between models to be a single relation on individuals 
that is a $q$-correlator for each $q$th order relation in the model. The 
resulting definition is more general, as well as correct.

The first time the notion of isomorphism for higher-type models is 
considered in print appears to be in \cite{Tarski1935a}.  However, 
it is not defined in nearly as much detail there as in Carnap's
\textit{Untersuchungen}.  The same definition occurs also in
\citet{LindenbaumTarski1935}, and there it coincides with
Carnap's.\footnote{It coincides, that is, except for a difference
Tarski stresses: his definition requires the isomorphism to be 
a bijection of the entire domain of individuals, not of just of 
individuals ``used'' in the model.}  The similarity of Carnap's 
and Tarski's approaches to isomorphism does not stop there.
In fact, Tarski's project of a general metamathematics of deductive
systems in the 1930s is similar to Carnap's in many other respects
as well. Tarski carries out his general investigations of ``deductive 
systems'', i.e., axiom systems, in a metatheory similar to Carnap's: 
the simple theory of types.  An axiom system for Tarski, as for Carnap, 
is a propositional function in the sense discussed above.

Indeed, as late as 1940 Tarski writes the following:
\begin{quote}
Let us consider a system of non-logical sentences and let, for
instance ``$C_1$'', ``$C_2$'', \dots, ``$C_n$'' be all the non‐logical
constants which occur. If we replace these constants by variables
``$X_1$'', ``$X_2$'', \dots, ``$X_n$'' our sentences are transformed
into sentential functions with $n$ free variables and we can say that
these functions express certain relations between $n$ objects or
certain conditions to be fulfilled by $n$ objects. Now we call a
system of $n$ objects $O_1$, $O_2$, \dots, $O_n$ a model of the
considered system of sentences if these objects really fulfill all
conditions expressed in the obtained sentential functions. [\dots] We
now say that a given sentence is a logical consequence of the system
of sentences if every model of the system is likewise a model of this
sentence. \citep{Tarski1940}
\end{quote}
This is the very same definition that Tarski gave in his 
seminal paper on logical consequence five years earlier 
\citeyearpar{Tarski1936}.\footnote{It bears emphasizing here that, 
while in that paper Tarski alluded to G\"odelian arithmetization of 
syntax, he did so only in order to explain how derivability from the
axioms of a formal system according to given rules can be 
formalized in metalogic; but this did not play a role in his definition 
of consequence.  In particular, Tarski did not give a uniform
definition of ``follows from'' for an object language in a formal
metalanguage which can refer to and quantify over all models and
expressions, in contrast to his general definition of truth.  Of course,
Tarski's discussion of consequence goes much beyond Carnap's and is
more nuanced in other respects, e.g., in its considerations of the distinction 
between derivability and consequence, in terms of which constants should 
count as logical, etc. He was able to accomplish that because G\"odel had 
in the meantime identified the gap between what follows from and what 
can be derived from deductive theories.  Tarski was also prompted to 
do it by Carnap's \emph{Logical Syntax}.}  In his 1940 talk Tarski
goes on to define categoricity in the same way in which Carnap defined
monomorphicity:
\begin{quote}
In order to obtain the second variant of the concept of categoricity,
we shall confine ourselves for the sake of simplicity, to such a case
in which the considered system of sentences is finite and contains
only one non‐logical constant, say ``$C$''. Let ``$P(C)$'' represent
the logical product of all these sentences. ``$C$'' can denote, for
instance, a class of individuals, or a relation between individuals,
or a class of such classes or relations etc. [\dots] We can now
correlate with the semantical sentence stating that our system of
sentences is semantically categorical an equivalent sentence
formulated in the language of the deductive theory itself. This is the
following sentence [\dots] in symbols:
\[
\forall X\forall Y(P(X) \land P(Y) \to X \sim Y)
\]
Now we say that the considered system of sentences is categorical with
respect to its logical basis \dots if the sentence formulated above is
logically valid.\footnote{Tarski uses $X \sim Y$ as an abbreviation
for the the isomorphism relation.}
\end{quote}
By this time---due to G\"odel's seminal work, but also Carnap's
contributions---the need for a proper way of talking about expressions 
of the axiom system in a metalanguage, in which general questions about
axiomatic systems are pursued (in particular, quantification over such
expressions to properly express derivability, consistency, and
completeness), was clear, and the proper solution (arithmetization) was
properly understood as well.  At the same time, what is evident from 
Tarski's writings on the subject into the 1940s is that a mature approach 
of model theory was not yet available; semantics was done very 
much in the way Carnap had proposed.

\section{a- and k-concepts}
\label{akConcepts}

As we saw earlier, Carnap's attempt to justify his identification of provability 
with (the derivability in type theory of) his statement expressing consequence 
was unsuccessful.  The problem here extends to his notion of consistency.  
Carnap does not properly locate the relevant definitions in the
metalanguage (so as to use an inductive, syntactic notion of proof).
Hence his proof of ``completeness'' does not establish the connection 
we expect today from such a theorem, namely one between proof 
theory and model theory.  At best, that proof can be seen as an
attempt to show that the existence of a model for an axiom system is
equivalent, in type theory, to the absence of a contradiction that is
a consequence of the axioms.  Carnap's approach also remained 
an unsuccessful attempt due to the problematic nature of the
type-theoretic quantification over propositional functions required
for formulating ``there is a contradiction which follows from the
axiom system''.

Having said that, Carnap was well aware of the difference between 
there being such a contradiction and actually having found
one.  This difference loomed large in the foundational debate between
Hilbertian ``formalists'' and Brouwerian intuitionists at the time.
Carnap attempted to capture it, and to do justice to it within his 
general axiomatics, by distinguishing ``a-concepts'' from ``k-concepts.''  
The latter are ``constructive [\emph{konstruktiv}]'' in the sense that
they require either a witness or a procedure to identify a witness; the 
former are ``absolutist'' in the sense that they only require the truth 
of an existential claim but not an actual witness. Carnap's a-concepts
can be inter-related within type theory by considering corresponding
conditional sentences.  (It can be proved that they are related
in certain ways by proving those sentences type-theoretically.)  For
instance, he is able to prove the equivalence of ``a-consistent''
and ``a-satisfied'' along such lines: the sentences $\exists \mR\,
f\mR$ and $\lnot \exists h\forall \mR(f\mR \to (h\mR \land \lnot
h\mR))$ are provably equivalent in type theory.\footnote{We leave
aside here, once again, the problems arising from taking $\exists h$ 
to correctly express the existence of a formula (``propositional
function'' in that sense).}

k-Concepts can also be so related, at least in some cases.  For 
that purpose, Carnap has to pick what he calls a ``criterion'' for a
particular k-concept to apply.  In the case of concepts involving
existential claims, Carnap chooses as criterion being able to exhibit
witnesses.  For instance, an axiom system is k-satisfied if a model
can be exhibited, and it is k-inconsistent if a propositional function
$h$ can be given that is provably formal and for which both $\forall
\mR(f\mR \to h\mR)$ and $\forall \mR(f\mR \to \lnot h\mR)$ are
provable.  When the k-concept involves a general claim, Carnap
requires that the formula expressing it is provable in type theory.
Thus, an axiom system is k-empty if $\lnot\exists \mR\, f\mR$ is
provable, and it is k-consistent if $\lnot\exists h\forall \mR(f\mR
\to (h\mR \land \lnot h\mR))$ is provable.  It follows---and Carnap
acknowledges this explicitly---that even when a partition 
into a-concepts is exhaustive (e.g., every axiom system is either
a-satisfied or a-empty), the corresponding k-concepts do not
necessarily form an exhaustive disjunction: e.g., an axiom system is
neither k-satisfied nor k-inconsistent if we neither have a model nor
a proof of a contradiction.

When Carnap considers a- and k-concepts in parallel, he usually
compares how properties of axiomatic systems can be defined
classically and how they can be defined intuitionistically, as well as
which results can be established classically and which
intuitionistically.\footnote{Behmann suggested that Carnap remain
  neutral between the absolutist and constructivist approaches
  (Behmann to Carnap, March 6, 1929, Staatsbibliothek zu Berlin,
  Nachla\ss{} 335 (Behmann), K.~1 I~10).  The first part (pp. 1--34)
  of Carnap's typescript RC 080--34--02 apparently is a revised
  version dating from after Behmann's letter: it does not bear
  Behmann's marginal annotations which are present on the rest of the
  typescript, and the version to which Behmann responded had \S10
  covering pp. 31--35, while in the surviving typescript \S10 covers
  pp.~30--35.}  In these comparisons, he comes close to making a
distinction that we now consider natural and important for metalogic:
that between truth and provability.  However, he does not quite get
there.  Even though, e.g., the definition of k-inconsistency involves
provability of a contradiction, it requires more than our notion of
inconsistency: in line with an intuitionistic (constructive)
interpretation, the contradiction must actually be exhibited.  More
importantly, the k-concepts corresponding to complements of a-concepts
are not the complements of the corresponding k-concepts.  This is so
because the constructive interpretations of properties of axiom
systems involving general claims require the provability in type
theory of the sentences expressing that an axiom system has those
properties (see above).  For example, while our notion of consistency
only requires the absence of a proof of a contradiction, Carnap's
k-consistency requires a proof in type theory of this absence.

Carnap realized that k-consistency does not imply
k-satisfaction, i.e., provability of 
\[\lnot\exists h\forall \mR(f\mR \to (h\mR \land \lnot h\mR))\] 
does not imply that a model can be explicitly given.  He could also 
have considered the question of whether k-consistency implies 
a-satisfiability, or more generally, the relationship between 
k-concepts and corresponding a-concepts.  He did not do so, 
perhaps because he assumed considerations involving a-concepts
and k-concepts, while parallel, to be independent, or maybe because it is
obvious that a k-concept implies the corresponding a-concept but never
vice versa.  From today's perspective, the interesting cases are
those that involve the relations between the mere failure of a k-concept 
to apply and the opposite k- or a-concept.  For example, completeness 
(in our sense) is the question of whether an axiom system that is
not k-inconsistent is k-satisfied or at least a-satisfied. These,
however, are simply not natural questions to ask in Carnap's
framework, since they require a clearer distinction between proving
things from and proving things about axiom systems in the
\emph{Grunddisziplin}.

Nevertheless, Carnap's discussion of k-concepts is 
illuminating with respect to his notion of ``decidability
[\emph{Entscheidungsdefinitheit}]''.  This notion is reminiscent 
of our notion of semantic completeness for theories, as we saw above, 
and his corresponding k-notion is related to both decidability and 
syntactic completeness of theories.  Carnap defines an axiom 
system~$f$ to be a-decidable as follows:
\[
\exists \mR \land \forall g(\For(g) \to \forall \mR(f\mR \to g\mR) 
\lor \forall \mR(f\mR \to \lnot g\mR))
\]
This is provably equivalent, within type theory, to $f$ being
non-forkable.  The intent of this definition is that every formal
statement is either a consequence of the axioms or its negation is.
It differs from our notion of semantic completeness in that the
$\forall g$ quantifier is problematic as a formalization of ``for
all formulas'' (see above), also in that $h$ is allowed to
be any propositional function of type theory (not just, say, first
order formulas in the language of $f$).  Whereas in other cases of
general properties Carnap adopts provability in type theory as the
criterion for the k-concept, here he defines it thus: ``If a model of
$f$ can be exhibited, and a procedure can be given, according to
which, for every given formal $g$ (with the same variables), either a
proof of $f \to g$ or a proof of $f \to \lnot g$ can be carried out in
finitely many steps.''  The resulting notion is, as Carnap points out,
not coextensive with k-monomorphicity or k-nonforkability.  He also
claims that this is the notion that is usually meant when other 
authors at the time write about decidability.  Moreover, in contrast to
monomorphicity, no axiom systems were known to be k-decidable 
to him, nor was it clear whether any axiom systems ever would be 
shown to be k-decidable or not.

This latter claim may seem puzzling.  After all, Langford's proof of 
the decidability of theories of order (from \citeyear{Langford1927a}) 
was already published at the time .  It is possible that Carnap 
simply did not know about Langford's work,\footnote{Carnap 
does not mention Langford's work himself.  It is mentioned by
\citet{Fraenkel1928}, although not in the context of completeness 
or decidability of axiom systems.}  But even if he did, he would not
have recognized it as providing a decision procedure for establishing
k-decidability for axiom systems of order, since Langford 
only showed that the axioms settle the truth values of \emph{first-order}
propositional functions involving the same variables as the axiom
system.  We now take it for granted that this is the question to be
settled; but at the time Carnap was writing, the almost universal
restriction to first-order languages that we typically assume was 
still some way off.  Certainly for Carnap, but also for most of his 
contemporaries, logic was still higher-order, and the consequences 
of the axiom systems of interest included their higher-order
consequences.\footnote{Note that many of the axiom systems in
\emph{Abriss} involve higher-order axioms, even primitives or
important defined concepts that are higher order. For instance, the
axiom system for geometry in \S34 has a second-order primitive
\textsf{ger}, for the class of lines, which are themselves classes of
points.}

Another reason for why the claim above may seem puzzling is the 
following: Carnap himself stresses that some axiom systems are known 
to be k-monomorphic, i.e., in type theory provably monomorphic, and 
hence, provably non-forkable by any formal sentence (not just first-order
sentences).  Consequently, a simple decision procedure would be to 
conduct an exhaustive search of all proofs in type theory until a proof 
of either $f \to g$ or $f \to \lnot g$ is found.  Although this approach 
seems obvious to us today, it was far from obvious at the time.  For
instance, \citet[p.~366]{Bernays1930} explicitly considers this
approach and rejects it in the following passage:
\begin{quote}
Notice that this requirement of being deductively closed [i.e., every
statement that can be formulated within the framework of the theory
is either provable or refutable] does not go as far as the
requirement that every question of the theory be \emph{decidable}.
The latter says that there should be a procedure for deciding for any
arbitrarily given pair of contradictory claims belonging to the theory
which of the two is provable (``correct'').
\end{quote}

The much more inclusive range of consequences of axioms that
Carnap considers (in full type theory) also plays an important role in
his discussion of the relationship between the k-decidability of axiom
systems and the logical decision problem in sections 3.7 and~3.8 of
\emph{Untersuchungen}.  The latter is the question of whether there is
a procedure which decides for each sentence of logic (i.e., type
theory) after finitely many steps whether the sentence is provable or
not.  As is well known, Hilbert considered this problem one of the
main open problems of logic \citep{Hilbert1928,HilbertAckermann1928}
and he was optimistic that a positive solution was forthcoming.  This
optimism was not universally shared, however; Carnap cites
\citet{Weyl1925} to that effect.  If the logical decision problem was
solved positively, Carnap argues in section~3.7 of
\emph{Untersuchungen},\footnote{\S29 of RC 080--34--03, p.~98ff.}
every monomorphic axiom system would be k-decidable.  In that event,
k-decidability would be coextensive with monomorphicity, and thus, 
according to Carnap, not independently useful, since it does not 
result in a new division of axiom systems into two classes.

As long as this problem is not solved positively, i.e., as long as
there are undecided logical sentences, no axiom system is k-decidable.
This is so because the undecided logical sentence~$r$ yields, for each
axiom system $f$, a propositional function $g$ which likewise is
undecided.  (Recall here that being k-decidable depends not only on the
existence of a decision procedure, but on such a decision procedure
being known---in keeping with the intuitionistic interpretation of
k-decidability.)  Take any formula~$h$ which is a consequence of~$f$,
such as $f$ itself if $f$ is formal, or a tautology involving only the
variables in $f$.  Then $g\mR = h\mR \land r$ follows from $f$ if $r$
holds, and $\lnot h$ follows from $f$ iff $r$ does not hold.  In fact,
type theory proves
\[
\exists \mR\, f\mR \to [\forall \mR(f\mR \to h\mR) \leftrightarrow r]
\]  
Thus the exhibition of a model for $f$ together with a decision
procedure that yields a proof either of $\forall \mR(f\mR \to g\mR)$
or of $\forall \mR(f\mR \to \lnot g\mR)$ in type theory would yield a
proof of $r$ or $\lnot r$ in type theory.

As with other issues tackled by Carnap in the \emph{Untersuchungen}, 
his considerations about decidability, as well as his distinction between 
absolutist and constructive versions of metalogical notions, remained 
fundamentally flawed.  They do, nevertheless, reveal significant insight 
into the subtleties involved, and they constitute a valiant effort to 
address them with the tools Carnap had available.

\section{Conclusion}
\label{Summary}

In light of our discussion in the previous sections, we are left with 
a mixed assessment of Carnap's early metatheoretic work.  Earlier 
criticisms of it (e.g. by Coffa and Hintikka) did miss their mark in various 
ways, since they were too coarse-grained and dismissive.  Nevertheless, 
Carnap's approach was limited in several important respects, as our 
discussion also made evident.  To assess Carnap's place in the history 
of logic properly, an account that balances both sides is needed.

On the positive side, Carnap addressed and made progress with 
a considerable number of problems in \emph{Untersuchungen} and 
\emph{Abriss}. First of all, notions such as consequence, isomorphism, 
categoricity, and completeness were not clearly defined at the time  
and he made steps towards making them precise.  Moreover, he did so 
by using the best tools available to him, namely simple (or deramified) 
type theory.  Second, Carnap began to work out ideas appropriate to 
the model-theoretic character of these notions.  In particular, he was 
aware of the kind of domain variation underlying them.  His explication 
of that idea looks unsatisfactory from today's perspective.  But again, 
it was arguably the most promising route to take given the state of the
development of logic at the time.  For these reasons, it seems fair to say that 
Carnap anticipated the kind of formal reasoning about models with 
which we are now familiar in substantive ways; he thus provided an 
early, noteworthy form of model theory within a type-theoretic setting.  
Third, he was sensitive to open problems in the foundations of 
mathematics at the time, e.g., the distinction between ``absolute'' 
and ``constructive'' concepts, and he tried to give a precise 
characterization of them as well.  Fourth, in his discussion of 
the different completeness properties of axiom systems, Carnap 
asked fruitful questions and tried to prove important theorems, 
including the \emph{Gabelbarkeitssatz}.

On the negative side, there are also several aspects in which Carnap's
attempts remained unsuccessful, even by his own lights.  He discussed
provability as a meta-theoretic notion defined ``outside'' the
framework of type theory, but his erroneous proof that provability and
consequence coincide prevented him from fully grasping the problems
posed by decidability.  The axiom systems he considered, as well
as his model-theoretic definitions of consequence, satisfaction, etc.,
are specified in the object language.  As we saw, this is not outright
impossible, and with Tarski he was in good company in using this 
approach.  Yet it does have implications of which Carnap was not
aware in 1928--29 and which limit the importance of his results.  Second, 
some of the notions discussed in his \emph{Untersuchungen} manuscript 
are too general to be useful, as later developments made clear.  This 
concerns, for instance, his notion of formality, which is replaced
today by specific language restrictions.  A third problem Carnap's
approach faces is that some of the central concepts underlying his
general axiomatics project are left underspecified and, thus, largely
unexplored.  This concerns especially the nature of his ``basic
discipline'', including notions such as truth and provability, but also 
the relationship between absolute and constructive metatheoretic
concepts, and in particular, between a-forkability and k-decidability.

Several issues concerning Carnap's early metatheory remain in need of 
further exploration.  One of them is the influence Carnap's work on general
axiomatics, and his \emph{Untersuchungen} in particular, had on the
historical development of logic: How, more broadly and in more depth, 
do his contributions relate to work by Tarski, Hilbert, Bernays, and G\"odel,
and other central figures?\footnote{See \cite{SchiemerReck2013} for a
few steps in that direction, earlier also \cite{AwodeyCarus2001},
\cite{AwodeyReck2002}, \cite{Goldfarb2005}, \cite{Reck2007}, and
\cite{Reck2011}, among others.}  A second issue concerns Carnap's 
own intellectual development: Are there points of contact or continuity
between his early model theory, as discussed in the present paper, and, 
e.g., his ``official'' work on semantics from the 1940s?  Third, more
questions can be raised about the remaining significance of Carnap
early metatheory from a technical point of view.  This concerns the
\emph{Gabelbarkeitssatz}, in particular, including recent research on
the so-called ``Fraenkel-Carnap property'' for second-order theories
stimulated by it.  The latter has led to some interesting partial results 
already, i.e., to the identification of several mathematical theories that 
have this metatheoretic property.\footnote{See \cite{WeaverGeorge2005} 
and \cite{WeaverPenev2011} for recent work on the Fraenkel-Carnap 
property for second-order theories.  See \cite{AwodeyReck2002b} for more
general results related to Carnap's Gabelbarkeitssatz.}  Finally, Carnap's 
conjecture that (higher-order) theories are categorical if and only if 
they are semantically complete has still not been decided in general;  
it thus remains ``one of the leading open questions in higher-order 
axiomatics'' \citep[p.~84]{AwodeyReck2002}.

\begin{acknowledgements}
Georg Schiemer's research was supported by the Austrian Science Fund,
projects J3158--G17 and P--27718. Richard Zach's research was
supported by the Social Sciences and Humanities Research Council of
Canada.

We are grateful to A. W. Carus and an anonymous referee for their
detailed and thoughtful comments on earlier drafts of this paper. We
would also like to thank the audiences at the 2013 workshop
\emph{Carnap on Logic} at the Munich Center for Mathematical
Philosophy, the Minnesota Center for Philosophy of Science, the 2014
Spring Meeting of the ASL, and the 2014 Society for the Study of the
History of Analytic Philosophy Annual Meeting.  Richard Zach
acknowledges the generous support of the Calgary Institute for the
Humanities and the Department of Philosophy at McGill University.
\end{acknowledgements}

\bibliographystyle{spbasic}      
\bibliography{Carnap}

\end{document}